\documentclass[11pt]{article}
\usepackage{graphicx}
\usepackage{amssymb}
\usepackage{epstopdf}
\usepackage{amsmath}
\usepackage{amsfonts}

\DeclareSymbolFont{AMSb}{U}{msb}{m}{n}
\DeclareMathSymbol{\N}{\mathbin}{AMSb}{"4E}
\DeclareMathSymbol{\Z}{\mathbin}{AMSb}{"5A}
\DeclareMathSymbol{\R}{\mathbin}{AMSb}{"52}
\DeclareMathSymbol{\Q}{\mathbin}{AMSb}{"51}
\DeclareMathSymbol{\I}{\mathbin}{AMSb}{"49}
\DeclareMathSymbol{\C}{\mathbin}{AMSb}{"43}

\begin{document}

\addtolength{\textheight}{0 cm}
\addtolength{\hoffset}{0 cm}
\addtolength{\textwidth}{0 cm}
\addtolength{\voffset}{0 cm}

\setcounter{secnumdepth}{5}
\newcommand{\1}{I\hspace{-1.5 mm}I}
\newtheorem{proposition}{Proposition}[section]
\newtheorem{theorem}{Theorem}[section]
\newtheorem{lemma}[theorem]{Lemma}
\newtheorem{coro}[theorem]{Corollary}
\newtheorem{remark}[theorem]{Remark}
\newtheorem{ex}[theorem]{Example}
\newtheorem{claim}[theorem]{Claim}
\newtheorem{conj}[theorem]{Conjecture}
\newtheorem{definition}[theorem]{Definition}
\newtheorem{application}{Application}
\newtheorem{corollary}[theorem]{Corollary}

\def\LX{{\cal L}(X)}
\def\LY{{\cal L}(Y)}
\def\LH{{\cal L}(H)}
 \def\ASD{{\cal L}_{\rm AD}(X)}
 \def\ASDY{{\cal L}_{\rm AD}(Y)}
\def\ASDH{{\cal L}_{\rm AD}(H)}
 \def\ASDP{{\cal L}^{+}_{\rm AD}(X)}
  \def\ASDYP{{\cal L}^{+}_{\rm AD}(Y)}
   \def\ASDHP{{\cal L}^{+}_{\rm AD}(H)}
\def\CX{{\cal C}(X)}
\def\CY{{\cal C}(Y)}
\def\CH{{\cal C}(H)}
\def\PX{{\cal A}(X)}
\def\PY{{\cal A}(Y)}
\def\PH{{\cal A}(H)}
\def\phi{{\varphi}}
\def\AH{A^{2}_{H}}

\title{Stability Estimates for Coefficients of Magnetic Schr\"odinger Equation From Full and Partial Boundary Measurements}
\author{Leo Tzou\thanks{This research was supported by the Doctoral Post-Graduate Scholarship from the Natural Science and Engineering Research Council of Canada.}\\
\small Department of Mathematics\\
\small University of Washington\\
\small Seattle Washington 98105
}
\maketitle
\section*{Abstract} 
In this paper we establish a $log\ log$-type estimate which shows that in dimension $n\geq 3$ the magnetic field and the electric potential of the magnetic Schr\"odinger equation depends stably on the Dirichlet to Neumann (DN) map even when the boundary measurement is taken only on a subset that is slightly larger than half of the boundary $\partial\Omega$. Furthermore, we prove that in the case when the measurement is taken on all of $\partial\Omega$ one can establish a better estimate that is of $log$-type. The proofs involve the use of the complex geometric optics (CGO) solutions of the magnetic Schr\"odinger equation constructed in \cite{sun uhlmann} then follow a similar line of argument as in \cite{alessandrini}. In the partial data estimate we follow the general strategy of \cite{hw} by using the Carleman estimate established in \cite{FKSU} and a continuous dependence result for analytic continuation developed in \cite{vessella}.\\\\

\begin{section}{Introduction}
\hspace{0.5 cm}Throughout this article we assume that the dimension $n\geq3$. Let $\Omega\subset \R^n$ be an open bounded set with $C^\infty$ boundary, we are interested in the magnetic Schr\"odinger operator
\[H_{W,q} := D^2 + W\cdot D +D\cdot W + W^2 + q\]
with real vector valued magnetic potential $W\in W^{2,\infty}(\Omega,\R^n)$ and the bounded electric potential $q\in L^\infty(\Omega)$. As usual, $D := -i\nabla$.
For simplicity, we assume throughout that for all $(W,q)$ under consideration $0$ is not an eigenvalue of the operator $H_{W,q}: H^2(\Omega)\cap H^1_0(\Omega) \to L^2(\Omega)$.\\
Let $\nu$ be the unit outer normal. Under the present assumptions, the Dirichlet problem
\[H_{W,q}u = 0\ \ \ \ \ \ \ \ \ \ u\mid_{\partial\Omega} = f\]
has a unique solution in $H^1(\Omega)$, and we can introduce the Dirichlet to Neumann (DN) map $\Lambda_{W,q}: H^{1/2}(\partial\Omega)\to H^{-1/2}(\Omega)$ associated with the magnetic Schr\"odinger operator $H_{W,q}$ by
\[\Lambda_{W,q}: f\mapsto (\partial_\nu + iW\cdot\nu)u\mid_{\partial\Omega}\]
The inverse problem under consideration is to recover information about the magnetic and electric potential from the DN map measured on a subset of the boundary. 

\hspace{0.5 cm}In the absence of the magnetic potential, the identifiability problem was solved by \cite{sylvester uhlmann} for when the measurement is taken on the whole boundary. Recently, Kenig-Sj\"ostrand-Uhlmann in \cite{KSU} showed that the same result holds even if the measurement is taken on possibly a very small subset of the boundary. The issue of stability without the magnetic potential was first addressed by Alessandrini in \cite{alessandrini} for the full data problem and later by Heck-Wang \cite{hw} when the data is measured on a subset that is slightly larger than half of the boundary.

\hspace{0.5 cm}In the presence of a magnetic potential, it was noted in \cite{sun} that the DN map is gauge invariant. Namely, given any ${\cal P}\in C^2(\bar\Omega)$ with ${\cal P}\mid_{\partial\Omega}=0$, one has $\Lambda_{W+\nabla {\cal P},q} = \Lambda_{W,q}$. Therefore, the magnetic potential is not uniquely determined by $\Lambda_{W,q}$. However, as was shown in \cite{FKSU}, the magnetic field $dW$ and electric potential $q$ are uniquely determined even if the measurement is taken only on a small part of the boundary. Furthermore, recently in \cite{salo} a method was given for reconstructing the magnetic field and electric potential under some regularity assumptions on the magnetic potential.

\hspace{0.5 cm}Following the above identifiability and reconstruction results, it is natural to ask whether small perturbations in the DN map would lead to small changes in the $dW$ and $q$ determined by $\Lambda_{W,q}$. This paper establishes a $log\ log$-type stability estimate for $dW$ and $q$ in the case when the measurement is taken only on a chosen subset of $\partial\Omega$. In the process we will also show that if one has full data measurements, the result can be improved to a $log$-type estimate. As mentioned before, when the magnetic potential is absent, the full and partial data estimates are established in \cite{alessandrini} and \cite{hw} respectively by using complex geometric optics (CGO) solutions to the Schr\"odinger equation that approximate plane waves. We follow a similar strategy except that in the presence of the magnetic field we need to use a richer set of CGO solutions studied in \cite{sun uhlmann} and \cite{salo}.

\hspace{0.5 cm}This article is organized into three parts. In part I we need to prove some fine properties of the CGO solutions that were not considered in previous studies. This is because the existing theory in \cite{sun uhlmann} and \cite{salo} are sufficient for identifiability and reconstruction results but a slightly more refined understand is necessary for establishing stability. The difficult issue is the following. Given an $M >0$, $R>0$ $p>n$ we consider the family of compactly supported vector fields 
\[{\cal W}(M,R) := \{W\in W^{1,p}(\R^n;\R^n) \mid \|W\|_{W^{1,p}}\leq M, \|div W\|_{L^\infty}\leq M\  and\  supp(W)\subset B_R\}\] 
and the family of compactly supported electric potentials 
\[{\cal Q}(M,R) := \{q\in L^\infty\mid \|q\|_{L^\infty}\leq M, supp(q)\subset B_R\}\] 
In order to prove stability, we need to show that the CGO solutions to the equation $H_{W,q}u = 0$ has remainder decaying uniformly for all $W\in {\cal W}(M,R)$, $q\in {\cal Q}(M,R)$. More precisely,
\begin{proposition}
Let $\sigma_0,\theta$ be positive numbers satisfying $\sigma_0 + \theta < \frac{1}{4n+6}$. For all $M>0$ and $R>0$ there exists constants $C>0$, $h_0 >0$, and $\epsilon >0$ depending on dimension, $p$, $\sigma_0$, $\theta$, $M$ and $R$ only such that for all $(W,q)\in {\cal W}(M,R)\times{\cal Q}(M,R)$, $\zeta \in \C^n$ with $\zeta\cdot\zeta = 0$ and $|\zeta|\geq \frac{1}{h_0}$, there exists solutions to $H_{W,q}u = 0$ in $\R^n$ of the form
\[u(x,\zeta) = e^{i\zeta\cdot x}(e^{i\chi_{|\zeta|}\phi^\sharp} + r(x,\zeta)),\ \ \ \  \ \ \ \|r(\cdot,\zeta)\|_{H^t_\delta}\leq C |\zeta|^{t-\epsilon},\ \ \ \ \ \ t\in[0,2]\]
where $\phi^\sharp$ is defined by
\[\phi^\sharp(x) = \frac{-\sqrt{2}\zeta}{2\pi |\zeta|}\cdot\int_{\R^2}\frac{ W^\sharp(x-y_1\frac{\sqrt{2}Re(\zeta)}{|\zeta|} -y_2\frac{\sqrt{2}Im(\zeta)}{|\zeta|})}{y_1+iy_2}\]
Here $\chi_{|\zeta|} = \chi(x/|\zeta|^{\theta})$ with $\chi(x)$ being a smooth function supported in the unit ball and is $1$ near zero, and $W^\sharp$ is the convolution of $W$ with the mollifier defined by 
\[W^\sharp(x) = \int |\zeta|^{n\sigma_0}\chi(y|\zeta|^{\sigma_0}) W(x - y) dy\]
\end{proposition}

\hspace{0.5 cm}In part II, we use the CGO solutions constructed above to prove stability for the full data problem.  The computation for the magnetic field stability is analogous to the one given for the electric potential in \cite{alessandrini} combined with some ideas in \cite{salo}.
 
\begin{theorem}
Let $\Omega$ be a bounded open subset of $\R^n$ with smooth boundary. For all $M>0$, there exists a $C> 0$, $\epsilon>0$ such that the estimate
\[\|I_{\Omega}d(W_1 - W_2)\|_{H^{-1}(\R^n)}\leq C\lbrace \|\Lambda_{W_1,q_1}-\Lambda_{W_2,q_2}\|_{\frac{1}{2},\frac{-1}{2}}^{1/2} + |log \|\Lambda_{W_1,q_1}-\Lambda_{W_2,q_2}\|_{\frac{1}{2},\frac{-1}{2}}|^{-\epsilon}\rbrace\]

holds for all $W_1, W_2\in W^{1,\infty}(\Omega)$ and $q_1,q_2\in L^\infty(\Omega)$ satisfying $\|W_l\|_{ W^{1,\infty}(\Omega)}\leq M$, $\|q_l\|_{L^\infty(\Omega)}\leq M$ ($l = 1,2$) and $W_1 = W_2$ on $\partial\Omega$. Here $I_\Omega(x)$ is the indicator function of $\Omega$.
\end{theorem}
The estimate for the electric potentials, however, is slightly more involved. Unlike the identifiability results in \cite{sun uhlmann} and \cite{FKSU}, the above theorem does not make the first order terms in the magnetic Schr\"odinger equation vanish. Therefore, complications would arise when one tries to establish the estimate for the (lower order) electric potentials in the presence of the (higher order) magnetic fields. To remedy this difficulty, we first show by using the Hodge decomposition that the $d$ operator on differential forms is in some sense "bounded invertible" when restricted to the right subspaces. Then we will combine this fact with the estimate we have for $d(W_1-W_2)$ to obtain the estimate for the electric potentials.

\begin{theorem}
Let $\Omega$ be a bounded open subset of $\R^n$ with smooth boundary. For all $M>0$, there exists a $C> 0$, $\epsilon>0$ such that the estimate
\[\|q_1 - q_2\|_{H^{-1}(\R^n)}\leq C\lbrace \|\Lambda_{W_1,q_1}-\Lambda_{W_2,q_2}\|_{\frac{1}{2},\frac{-1}{2}}^{1/2} + |log \|\Lambda_{W_1,q_1}-\Lambda_{W_2,q_2}\|_{\frac{1}{2},\frac{-1}{2}}|^{-\epsilon}\rbrace\]

holds for all $W_1, W_2\in W^{2,\infty}(\Omega)$ and $q_1,q_2\in L^\infty(\Omega)$ satisfying $\|W_l\|_{ W^{2,\infty}(\Omega)}\leq M$, $\|q_l\|_{L^\infty(\Omega)}\leq M$ ($l = 1,2$) and $W_1\mid_D = W_2\mid_D$ on $\partial\Omega$. 

\end{theorem}
In part III we assume knowledge of the DN map on a subset of $\partial\Omega$ that is only slightly larger than half of the boundary and prove a stability result that is weaker than the ones above. To give a precise statement of the theorems would require more defintions and therefore they will be stated in the introduction section of part III. The proof follows the idea employed in \cite{hw} where one uses a Carleman estimate that is established in \cite{FKSU} to help suppress the missing piece of information and obtain an estimate for the difference of the Fourier transform on a wedge in phase space. To extend the estimate from the wedge to a ball we will use a continuous dependence result for analytic continuation developed in \cite{vessella}. After this is established, the stability result for the magnetic field and electric potential would follow by similar calculations done in part II.

\end{section}

\begin{section}*{PART I - Fine Properties of CGO Solution}
\begin{section}{Properties of Transport Equations}
To establish stability we will need to construct complex geometric optics (CGO) solutions to $H_{W_l,q_l}u_l = 0$ that are of the form
\[u_l = e^{i\zeta_l\cdot x}(e^{i\chi_{|\zeta|}\phi_l^\sharp} + r_l(x,|\zeta|))\]
where $\zeta_l\cdot\zeta_l = 0$ and $\phi_l^\sharp$ satisfies the transport equation 
\[-\mu\cdot\nabla\phi_l^\sharp = \mu_l\cdot W_l^\sharp\]
here $\mu_l\in\C^n$ is defined by $\zeta_l = \frac{|\zeta_l|}{\sqrt{2}}\mu_l$ and $W_l^\sharp$ is the convolution of $W_l$ with a mollifier. In this section we will collect some properties regarding how $\phi^\sharp$ depends on $W_l$ and the unit vector $\mu_l$. Throughout this article we will denote by $N_\mu^{-1}$ to be the inverse of the operator $\mu\cdot\nabla$. More precisely,
\[N_\mu^{-1} (f) := \frac{1}{2\pi}\int_{R^2}\frac{f(x - Re(\mu)y_1 - Im(\mu)y_2)}{y_1 + iy_2}dy_1 dy_2\]
for all $f\in L^{\infty}_c(\R^n)$. The general properties of this operator is summarized in the following two lemmas which we will state without proof. Interested reader can see \cite{salo}.
\begin{lemma}
\label{L-infinity bound}
Let $f\in W^{k,\infty}_0(\tilde\Omega)$ with $f = 0$ for $|x| > R$. Then $u = N^{-1}_\mu f \in W^{k,\infty}(\R^n)$ solves the equation $\mu\cdot\nabla u = f$ in $\R^n$ and satisfies for all multi-index $|\alpha|\leq k$
\[|\partial^\alpha u(x)|\leq C(\tilde\Omega)\|\partial^\alpha f\|_{L^\infty(\tilde\Omega)}(1 + |x_T|^2)^{-1/2}I_{B_R}(x_\perp)\]
where $x_T$ is the projection of $x$ to the plane $T = span\{Re(\mu),Im(\mu)\}$, $x_\perp = x - x_T$ and $I_{B_R}$ is the indicator function of the ball of radius $R$ around $0$.
\end{lemma}
Sometimes we need a version of lemma \ref{L-infinity bound} where $f$ and $\mu$ depend on a parameter. Let $V\subset\R^n$ be an open set and let $\gamma_j(\xi)\ (j=1,2)$ be a $C^\infty$ function of $\xi\in V$ which satisfy 
\[1-\epsilon \leq |\gamma_j(\xi)|\leq 1+\epsilon,\ \ \ |\gamma_1(\xi)\cdot\gamma_2(\xi)|\leq \epsilon\]
and also $|\partial^\alpha\gamma_j(\xi)|\leq M_1$ for $|\alpha|\geq 1$. 
\begin{lemma}
\label{local transport estimate}
Let $\epsilon >0$ be small enough and let $f(x,\xi)\in C^\infty(R^n\times V)$ satisfy $f(x,\xi) = 0 $ for $|x| > R$. Then the function
\[u(x,\xi) = \frac{1}{2\pi}\int_{\R^n}\frac{1}{y_1 + iy_2}f(x-y_1\gamma_1(\xi) - y_2\gamma_2(\xi),\xi)dy_1dy_2\]
is in $C^\infty(\R^n\times V)$ solves $(\gamma_1(\xi) + \gamma_2(\xi))\cdot\nabla_x u = f$ in $\R^n$ and satisfies
\[|\partial_x^\alpha\partial_\xi^\beta u(x,\xi)|\leq C_{\alpha,\beta,R,M_1}(\sum\limits_{|\gamma+\delta|\leq |\alpha +\beta|}\|\partial_x^\gamma\partial_\xi^\delta f\|_{L^\infty(\R^n\times V)})<x_T>^{|\beta| - 1}\chi_{B_M}(x_\perp)\]
where $x_T$ is the projection of $x$ to the plane $T = span\{\gamma_1(\xi),\gamma_2(\xi)\}$ and $x_\perp = x - x_T$
\end{lemma}
In proving stability we will be interested in the dependence of the $N_{\mu}^{-1}(-\mu\cdot W)$ on the parameter $\mu$. The next lemma states that the dependence is continuous provided that $W$ behaves reasonably well.
\begin{lemma} 
\label{uniform estimate}
Let $W\in C^{0,t}_c(\tilde\Omega)$ with $\|W\|_{C^{0,t}_c(\tilde\Omega)}\leq M$. For any $\mu_1,\mu_2\in S^{n-1} + iS^{n-1}$ such that $Re(\mu_l)\perp Im(\mu_l)$ for $l= 1,2$ we have the following estimate
\[\|N_{\mu_1}^{-1}(-\mu_1\cdot W) - N_{\mu_2}^{-1}(-\mu_2\cdot W)\|_{L^{\infty}(\tilde\Omega)}\leq C|\mu_1 - \mu_2|^t\]
and the constant depends only on the size of $\tilde\Omega$ and is uniform for all $\|W\|_{C^{0,t}_c(\tilde\Omega)}\leq M$.
\end{lemma}
\noindent{\bf Proof}\\
Since $W$ is supported in the bounded set $\tilde\Omega$ and $Im(\mu_l)\perp Re(\mu_l)$, we have that there exists an $R >0$ such that $W(x - Re(\mu_l)y_1 - Im(\mu_l)y_2) = 0$ for all $x\in\tilde\Omega$ and $|y_1e_1 + y_2e_2| \geq R$. Therefore we have for all $x\in\tilde\Omega$
\begin{eqnarray*}
&&|N_{\mu_1}^{-1}(-\mu_1\cdot W)(x) - N_{\mu_2}^{-1}(-\mu_2\cdot W(x))|\\
&&=\frac{1}{2\pi}|\int_{B_R}\frac{\mu_2\cdot W(x - Re(\mu_2)y_1 - Im(\mu_2)y_2)- \mu_1\cdot W(x - Re(\mu_1)y_1 - Im(\mu_1)y_2)}{y_1 + iy_2}dy_1 dy_2|\\ 
&&\leq C\int_{B_R}\frac{|\mu_2- \mu_2| \|W\|_{L^{\infty}}+ \|W\|_{C^{0,t}(\tilde\Omega)}(|Re(\mu_1 - \mu_2)y_1|+ |Im(\mu_1 - \mu_2)y_2|)^t }{y_1 + iy_2}dy_1 dy_2\\
&&\leq C R|\mu_1 - \mu_2|\|W\|_{L^{\infty}} + CR^2 \|W\|_{C^{0,t}(\tilde\Omega)}|\mu_2 - \mu_2|^t
\end{eqnarray*}
So we have $\|N_{\mu_1}^{-1}(-\mu_1\cdot W) - N_{\mu_2}^{-1}(-\mu_2\cdot W)\|_{L^{\infty}(\tilde\Omega)}\leq C|\mu_1 - \mu_2|^t$ as desired. $\square$\\\\
Combining lemma \ref{uniform estimate} and lemma \ref{L-infinity bound} we have the following corollary which will be useful later on.
\begin{corollary}
Under the same hypothesis as lemma \ref{uniform estimate} we have the following estimate on the exponential
\label{exponential estimate}
\[\|e^{iN_{\mu_1}^{-1}(-\mu_1\cdot W)} - e^{iN_{\mu_2}^{-1}(-\mu_2\cdot W)}\|_{L^\infty(\tilde\Omega)}\leq C(M,\tilde\Omega)|\mu_1 -\mu_2|^t\]
where $C>0$ depends only on the size of $\tilde\Omega$ and $M$.
\end{corollary}
\noindent{\bf Proof}\\
By lemma \ref{L-infinity bound} there exists an $M'$ such that $\|N_\mu^{-1}(-\mu\cdot W)\|_{L^\infty(\tilde\Omega)}\leq M'$ whenever $\|W\|_{L^\infty(\tilde\Omega)}\leq M$. Since the map $z\mapsto e^{iz}$ is Lipchitz in the closure of the ball $B_{M'}\subset \C$, we have that 
\[|e^{iN_{\mu_1}^{-1}(-\mu_1\cdot W)(x)} - e^{iN_{\mu_2}^{-1}(-\mu_2\cdot W)(x)}|\leq C(M')|N_{\mu_1}^{-1}(-\mu_1\cdot W)(x) - N_{\mu_2}^{-1}(-\mu_2\cdot W)(x)|\]
for all $x\in\tilde\Omega$. Now apply lemma \ref{uniform estimate} we have the desired estimate for the exponential.$\square$\\

Let $\eta_\epsilon$ be the standard mollifier and denote by $W^\sharp = \eta_\epsilon*W$. The next lemma tells us how well $e^{iN_\mu^{-1}(-\mu\cdot W^\sharp)}$ approximates $e^{iN_\mu^{-1}(-\mu\cdot W)}$.
\begin{lemma}
\label{approximation by mollifier}
Let $W \in C^{0,t}_c(\tilde\Omega;\R^n)$ with $\|W\|_{C^{0,t}_c(\tilde\Omega;\R^n)}\leq M$. We have the following estimate
\[\|e^{iN_\mu^{-1}(-\mu\cdot W^\sharp)} - e^{iN_\mu^{-1}(-\mu\cdot W)}\|_{L^{\infty}(\tilde\Omega)}\leq Ce^{2CM}M\epsilon^t \]
for all $\mu\in S^{n-1}+iS^{n-1}$ with orthonormal unit real and imaginary part.
\end{lemma}
\noindent{\bf Proof}\\
Pick $R > 0$ large enough such that for all $x\in\tilde\Omega$ $W(x - y_1Re(\mu) - y_2Im(\mu))= 0$ whenever $|y_1 e_1 + y_2 e_2| > R$. Then we have for all $x\in \tilde\Omega$,
\begin{eqnarray*}
|e^{iN_\mu^{-1}(-\mu\cdot W^\sharp)(x)} - e^{iN_\mu^{-1}(-\mu\cdot W)(x)}|&\leq& e^{2M}|N_\mu^{-1}(-\mu\cdot W^\sharp) - N_\mu^{-1}(-\mu\cdot W)|\\&\leq& e^{2M}\int_{B_R}\frac{|(W^\sharp - W)(x - y_1Re(\mu) - y_2Im(\mu))|}{|y_1 + iy_2|}\\&&\leq \epsilon^te^{2M}M\int_{B_R}\frac{1}{|y_1 + iy_2|}
\end{eqnarray*}
so the lemma is complete. $\square$

The following result on nonlinear Fourier transform was used by Salo \cite{salo} in reconstruction methods. We will repeat it here for convenience of the reader. Similar ideas appear in Sun \cite{sun} and Eskin-Ralston \cite{ralston}
\begin{lemma}
\label{nonlinear fourier}
Assume that $\gamma\perp\tilde\gamma\perp\xi$ with $\gamma,\tilde\gamma\in S^{n-1}$ and define $\mu = \gamma + i\tilde\gamma$. Define $\Phi(x) := N^{-1}_\mu(-\mu\cdot W)$ with vector field $W\in L^{\infty}_c$. Then we have the following identity for the nonlinear Fourier transform:
\[\int_{\R^{n}}\mu\cdot W e^{i\xi\cdot x}e^{i\Phi(x)}dx =\int_{\R^{n}}\mu\cdot W e^{i\xi\cdot x}dx\]
\end{lemma}
\noindent{\bf Proof}\\
Without loss of generality, we can assume that $\mu = e_1 + ie_2$ since the general case can be reduced to this case via an orthonormal linear transform. With this choice of $\mu$, $\xi = (0,0,\xi')$ with $\xi'\in\R^{n-2}$.
\[\int_{\R^n}W(x)\cdot\mu e^{i\xi\cdot x}e^{i\Phi(x)}dx = \int_{\R^n}(-(\partial_1 +i\partial_2)\Phi(x))e^{i\Phi}e^{i\xi'\cdot x'}dx = \int_{\R^{n-2}}e^{i\xi'\cdot x'}h(x')dx'\]
where
\begin{eqnarray*}
h(x') &=& i\int_{\R^2}(\partial_1 +i\partial_2)e^{i\Phi(x_1,x_2,x')}dx_1dx_2 =\lim_{R\to\infty}i\int\limits_{|x_1e_1 + x_2e_2|\leq R}(\partial_1 +i\partial_2)e^{i\Phi(x_1,x_2,x')}dx_1dx_2\\
&=&i\lim_{R\to\infty}\int\limits_{|x_1e_1 + x_2e_2|= R}e^{i\Phi(x_1,x_2,x')}(\nu_1 + i\nu_2)dS(x_1,x_2)
\end{eqnarray*}
When $|x_1e_1 + x_2e_2|$ gets large, $e^{i\Phi(x_1,x_2,x')} = 1 + i\Phi + O(|\Phi|^2) = 1 + i\Phi + O(|x_1e_1 + x_2e_2|^{-2})$ by lemma \ref{L-infinity bound} so
\begin{eqnarray*}i\lim_{R\to\infty}\int\limits_{|x_1e_1 + x_2e_2|= R}e^{i\Phi(x_1,x_2,x')}(\nu_1 + i\nu_2)dS &=& -\lim_{R\to\infty}\int\limits_{|x_1e_1 + x_2e_2|= R}\Phi(x_1,x_2,x')(\nu_1 + i\nu_2)dS\\
&=&-\int_{|x_1e_1 + x_2e_2|\leq R}(\partial_1 + i\partial_2)\Phi(x_1,x_2,x')dx_1dx_2\\
&=&\int_{|x_1e_1 + x_2e_2|\leq R}\mu\cdot W(x_1,x_2,x')dx_1dx_2
\end{eqnarray*}
and the proof is complete.$\square$
\end{section}

\begin{section}{Semiclassical pseudodifferential calculus}
The results which appear in \cite{salo}, \cite{sun}, and \cite{sun uhlmann} rely on solutions to $H_{W,q}u = 0$ that are of the form
\[u(x,\zeta) = e^{\zeta\cdot x}(e^{i\chi_{|\zeta|} \phi^\sharp} + r(x,\zeta))\]
where $\phi$ is the solution of some transport equation and $r(x,\zeta)$ satisfies
\[ \|r(\cdot,\zeta)\|_{H^t_\delta}\leq C(W,q)|\zeta|^{t-\epsilon},\ \ \ \ \ t\in[0,2]\]
The situation in establishing stability is more delicate, however, since one considers a family of magnetic and electric potentials satisfying certain a-priori estimates. Therefore, more care is needed if we wish to establish a stability estimate that is uniform for all magnetic potentials under consideration.
In particular, we need to ensure that the constant $C$ which appears in the above estimate for the remainder term $r(x,\zeta)$ is uniformly bounded for all $W$ satisfying our a-priori assumption. It is with this in mind that we develop some explicit estimates for semi-classical $\Psi DO$ in terms of its symbols. Most of the results in this section are well know but we will nevertheless include proofs or exact references for completeness. All integration of symbols against complex exponentials are understood to be oscillatory integrals (see \cite{raymond}). We begin with a fundamental result which gives a sharp estimate on the operator norm of the $\Psi DO$ by its symbol:
\begin{proposition}(Calderon-Vaillancourt)
\label{operator by symbol}
There exists a $k(n)\in\N$ depending on dimension of $\R^n$ only such that for all semiclassical symbols $a(x,\xi; h)\in S_\sigma^0$ with $0\leq\sigma\leq 1/2$, the following estimate on $\|Op_h(a)\|_{L^2(\R^n)\to L^2(\R^n)}$ holds for all $0<h\leq 1$
\[\|Op_h(a)\|_{L^2(\R^n)\to L^2(\R^n)}\leq (2\pi)^{-n} \|{\cal B}\|_{TR}\sum\limits_{|\alpha|,|\beta|\leq k(n)} p^{\sigma}_{\alpha,\beta}(a)\]
where $p_{\alpha,\beta}^{\sigma}$ is a semi-norm on $S^0_{\sigma}$ defined by $p^\sigma_{\alpha,\beta}(a) := \sup\limits_{x,\xi\in\R^n,0< h\leq 1}\{h^{\sigma(|\alpha|+|\beta|)}|\partial^\alpha\partial^\beta a(x,\xi;h)|\}$ and ${\cal B}$ is the $\Psi DO$ of trace class defined by the symbol ${\cal B}(x){\cal B}(\xi)$ with $\hat{\cal B}(\xi) = <\xi>^{-2k(n)}$ .
\end{proposition}
\noindent{\bf Proof}\\
The proposition can be reduced to proving estimates for classical $\Psi DO$ of symbol order zero. In particular, the classical result by Calderon (see for example p.10 Vol. II of \cite{taylor2}) states that for any classical symbol $a(x,\xi)$ of order zero we have the following estimate for the corresponding operator $A$: 
\begin{eqnarray}
\label{calderon}
\|A u\|_{L^2(\R^n)\to L^2(\R^n)}\leq (2\pi)^{-n} \|{\cal B}\|_{TR}\|a(x,\xi)\|_{C^{k}(\R^n_x,\R^n_\xi)}
\end{eqnarray}
Now for all $a(x,\xi;h)\in S^0_\sigma$ define for each fixed $h>0$ the (classical) symbol $a_h(x,\xi) := a(\sqrt{h}x, \sqrt{h}\xi;h)$ of order zero. Observe that we have the following relationship between the semi-classical quantization $Op_h(a)$ and the classical quantization of $a_h(x,\xi)$:
\[(Op_h(a)u)(x) = \sqrt{h}^{-n}A_hu_h(x)\]
where $u_h(x)$ is defined by its Fourier transform $\hat u_h(\xi) := \hat u(\frac{\xi}{\sqrt{h}})$ and $A_h$ is the operator associated to the (classical) symbol $a_h(x,\xi)$. The proposition then follows after some simple calculation by applying (\ref{calderon}) to $A_h$ and use the fact that $\sigma \leq 1/2$.
$\square$\\
We will hence forth denote by $k(n)$ to be the smallest integer for which proposition \ref{operator by symbol} holds. Now we will derive a result regarding the composition of semi-classical $\Psi DO$. It is well known fact that if $a,b\in S^0_\sigma$ then $Op_h(b)Op_h(b) = Op(ab) + h^{1-2\sigma}Op_hS_\sigma^0$. However, we need an explicit estimate of the semi-norms of the remainder by the semi-norms of $a,b$. The next lemma establishes this but is only interesting for $0<\sigma <\frac{1}{4n + 6}$.
\begin{lemma}
\label{multiplication remainder}
Let $a(x,\xi;h), b(x,\xi; h)\in S^0_\sigma$ be semiclassical symbols. We have the identity 
\[Op_h(a) Op_h(b) = Op_h (c)\] 
where $c\in S^0_\sigma$ satisfies 
\[c(x,\xi; h) = a(x,\xi; h)b(x,\xi; h) + h^{1-\sigma(4n+6)}m(x,\xi; h)\]
with $m(x,\xi;h)\in S^0_\sigma$. Furthermore, the semi-norms of $m(x,\xi;h)$ satisfies
\[p_{\alpha,\beta}^\sigma(m) \leq C_{\alpha,\beta} (\sum\limits_{|\alpha'|,|\beta'|\leq |\alpha|+|\beta| + 2n+3}p_{\alpha',\beta'}^\sigma(a) + \sum\limits_{|\alpha'|,|\beta'|\leq |\alpha|+|\beta| + 2n+3}p_{\alpha',\beta'}^\sigma(b))\]
The constant $C_{\alpha,\beta}>0$ depends only on $\alpha,\beta$ and the dimension.
\end{lemma}
\noindent{\bf Proof}\\
Simple calculation yields that
\[Op_h(a)Op_h(b) u(x) = \frac{1}{(2\pi)^n}\int_{\R^n}c(x,h\xi;h) e^{ix\cdot\xi}\hat u(\xi)d\xi\]
where
\[c(x,\xi;h) := \frac{1}{(2\pi h)^n}\int_{\R^n}\int_{\R^n}e^{-i\frac{<x-y,\xi-\eta>}{h}}a_h(x,\eta)b_h(y,\xi) d\eta dy\]
For the sake of clarity, in this calculation we denote by $a_h(x,\xi) := a(x,\xi;h)$. Make change of variable $y' = y-x$ 
\[c(x,\xi;h) := \frac{1}{(2\pi h)^n}\int_{\R^n}\int_{\R^n}e^{i\frac{<y',\xi-\eta>}{h}}a_h(x,\eta)b_h(y'+x,\xi) d\eta dy\]
Taylor formula gives $b_h(y'+x,\xi) = b_h(x,\xi) + \sum\limits_{|\alpha|=1}y'^\alpha \int_0^1(\partial^\alpha_x b_h)(x + \theta y',\xi)d\theta$. Substitute this into the above equation we get
\[c(x,\xi;h) := \frac{1}{(2\pi h)^n}\int_{\R^n}\int_{\R^n}e^{i\frac{<y',\xi-\eta>}{h}}a_h(x,\eta)(b_h(x,\xi) + \sum\limits_{|\alpha|=1}y'^\alpha \int_0^1(\partial^\alpha_x b_h)(x + \theta y',\xi)d\theta) d\eta dy\]
The first term in this can be computed explicitly by using Fourier and inverse Fourier formula.
\begin{eqnarray}
\label{main term}
\frac{1}{(2\pi h)^n}\int_{\R^n}\int_{\R^n}e^{i\frac{<y',\xi-\eta>}{h}}a_h(x,\eta)(b_h(x,\xi)d\eta dy = b(x,\xi;h) a(x,\xi;h)
\end{eqnarray}
The next term will be the remainder which we will compute as explicitly as possible
\begin{eqnarray*}
\frac{1}{(2\pi h)^n}\int_{\R^n}\int_{\R^n}e^{i\frac{<y',\xi-\eta>}{h}}a_h(x,\eta)\sum\limits_{|\alpha|=1}y'^\alpha \int_0^1(\partial^\alpha_x b_h)(x + \theta y',\xi)d\theta d\eta dy'
&=& h^{1-(4n+6)\sigma}m(x,\xi; h)
\end{eqnarray*}
where $m(x,\xi;h)$ is given by the formula
\begin{eqnarray*}
&&m(x,\xi;h) := \frac{-i h^{(4n+6)\sigma}}{(2\pi)^n}\times\\
&&\ \ \ \ \ \ \ \sum\limits_{|\alpha|=1}\int_{\R^{2n}} \frac{e^{i<y,\eta>}((I-\Delta_\eta)^{n+1}\partial_\eta^\alpha a_h)(x,\eta+\xi)(I-\Delta_y)^{n+1}\int_0^1(\partial^\alpha_x b_h)(x + \theta h y,\xi)d\theta}{(1+|\eta|^2)^{n+1}(1 + |y|^2)^{n+1}} d\eta dy
\end{eqnarray*}
So $c(x,\xi;h)$ can be written as $a(x,\xi;h) b(x,\xi;h) + h^{1-(2n+4)\sigma}m(x,\xi;h)$. It remains to check that $m(x,\xi;h)\in S^{0}_\sigma$ and satisfies the seminorm estimates stated in the lemma. Observe that for all multi-indices $\beta,\gamma$,
\begin{eqnarray}
\label{seminorm}
\sup\limits_{\eta,x,\xi,h>0} | h^{(|\beta|+|\gamma|+2n+3)\sigma}(I-\Delta_\eta)^{n+1}\partial_\xi^\gamma\partial_x^\beta\partial_\eta^\alpha a_h(x,\eta+\xi;h)| \leq \sum\limits_{\alpha', \beta'\leq |\beta|+|\gamma|+2n+3} p^\sigma_{\alpha',\beta'}(a)
\end{eqnarray}
and the same holds for $b(x,\xi;h)$. The term involving the Laplacian in $y$ is
\[(I-\Delta_y)^{n+1}\big\{\frac{\int_0^1(\partial^\alpha_x b_h)(x + \theta h y,\xi)d\theta}{(1 + |y|^2)^{n+1}}\big\}\]
Taking derivatives of $\frac{1}{(1+|y|^2)^{n+1}}$ results only in more decay, therefore
\[h^{(|\beta|+|\gamma|+2n+3)\sigma}(I-\Delta_y)^{n+1}\partial_\xi^\gamma\partial_x^\beta\big\{\frac{\int_0^1(\partial^\alpha_x b_h)(x + \theta h y,\xi)d\theta}{(1 + |y|^2)^{n+1}}\big\} \leq C_n\frac{ \sum\limits_{\alpha', \beta'\leq |\beta|+|\gamma|+2n+3} p_{\alpha', \beta'}(b)}{(1+|y|^2)^{n+1}}\]
where $C_n$ depends only on the dimension $n$.
Combining the above inequality and (\ref{seminorm}) in addition to the fact that $\frac{1}{(1+|y|^2)^{n+1}}$ is integrable in $\R^n$, we obtain directly from the definition of $m(x,\xi;h)$ that for all multi-indices $\beta,\gamma$,
\begin{eqnarray*}
h^{\sigma(|\gamma|+|\beta|)}\partial_\xi^\gamma\partial_x^\beta m(x,\xi;h) &\leq& C_{n,\alpha,\beta}(\int_{\R^n}\frac{1}{(1+|y|^2)^{n+1}}dy)^2\times\\&&\ \ \ \ \ \ \ \ \ \ \ (\sum\limits_{\alpha', \beta'\leq |\beta|+|\gamma|+2n+3} p^\sigma_{\alpha',\beta'}(a)+ \sum\limits_{\alpha', \beta'\leq |\beta|+|\gamma|+2n+3} p^\sigma_{\alpha',\beta'}(b))
\end{eqnarray*}
holds for all $x,\xi\in\R^n$ and $0< h \leq1$ and the constant $C_{n,\gamma,\beta}$ depends only on dimension and multi-indices $\gamma,\beta$.
$\square$\\
Given $a\in S^0_\sigma$ one can define a formal adjoint to $Op_h(a)$ in the usual way. It turns out that $Op_h(a)^*$ is also a semi-classical $\Psi DO$ with symbol denoted by $a^*\in S_\sigma^0$. Finer properties of this symbol is developed in the next lemma 
\begin{lemma}
\label{adjoint symbol estimate}
Let $a\in S^0_\sigma$ for $\sigma \leq 1/2$. Then $Op_h(a)^* = Op_h(a^*)$ with $a^*\in S^0_\sigma$. Furthermore, semi-norms of $a^*(x,\xi,h)$ can be bounded by semi-norms of $a(x,\xi,h)$ in the sense that for all multi-indices $\alpha, \beta$
\[p_{\alpha,\beta}^\sigma (a^*) \leq C_{\alpha,\beta}\sum\limits_{|\alpha'|,|\beta'|\leq |\alpha| + |\beta| + 2n+2}p_{\alpha',\beta'}(a)\]
with $C_{\alpha,\beta}$ depending only on the multi-indices and dimension.
\end{lemma}
\noindent{\bf Proof}\\
For each fixed $h>0$, set $a_h(x,\xi) := a(x,h\xi;h)$ and define the (classical) pseudodifferential operator 
\[A_h u := \int e^{i<x,\xi>}a_h(x,\xi) \hat u(\xi)d\xi = Op_h(a) u\]
With this notation we have for the formal adjoint $Op_h(a)^* = A_h^*$. By the result in classical pseudodifferential operator, $A_h^*$ can be written as
\[A_h^* u = \int e^{i<x,\xi>}a^*(x,h\xi;h) \hat u(\xi)d\xi\]
provided we set
\[a^*(x,\eta;h ) = \frac{1}{(2\pi h)^n}\int\int e^{i\frac{<z,\xi>}{h}} a(x+z, \xi - \eta;h)d\xi dz\]
It now remains to show that $a^*(x,\eta;h)\in S^0_\sigma$ and that its semi-norms satisfy the desired estimates. Taking $(I - h\Delta)^{n+1}$ of the exponential and integrate by parts as in the previous lemma we obtain
\[a^*(x,\eta;h ) = \frac{1}{(2\pi h)^n}\int\int \frac{e^{i\frac{<z,\xi>}{h}}}{(1 +\frac{1}{h}|\xi|^2)^{n+1}}(I - h\Delta_z)^{n+1}\frac{(I - h\Delta_\xi)^{n+1}a(x+z, \xi - \eta;h)}{(1 +\frac{1}{h}|z|^2)^{n+1}}d\xi dz\]
Make the change of variable $z' = \frac{z}{\sqrt{h}}$, $\xi' = \frac{\xi}{\sqrt{h}}$ and observing that $(I - h\Delta_z) = (I - \Delta_{z'})$ we get 
\[a^*(x,\eta;h ) = \frac{1}{(2\pi )^n}\int\int \frac{e^{i<z,\xi>}}{(1 +|\xi|^2)^{n+1}}(I - \Delta_{z'})^{n+1}\frac{(I - \Delta_{\xi'})^{n+1}a(x+\sqrt{h}z', \sqrt{h}\xi' - \eta;h)}{(1 +|z|^2)^{n+1}}d\xi' dz'\]
The semi-norm estimates now follows by similar arguments used in the previous lemma.
$\square$\\
Sometimes it is useful to conjugate $Op_h(a)$ with $<h D>^s$ to produce another semi-classical $\Psi DO$ of class $Op_h(S_\sigma^0)$.
\begin{lemma}
\label{conjugate}
Let $a\in S_\sigma^0$. For $|s| \leq 3$ we have $<hD>^{-s} Op_h(a) <hD>^s = Op_h(b)$ with $b\in S_\sigma^0$ and semi-norms of $b$ are bounded by semi-norms of $a$ in the usual sense:
\[p_{\alpha,\beta}^\sigma(b) \leq C_{\alpha,\beta}\sum\limits_{|\alpha'|,|\beta'|\leq n+3 + |\alpha|+|\beta|} p_{\alpha',\beta'}^\sigma(a)\] 
\end{lemma}
\noindent{\bf Proof}
Simple computation shows that we can write
\[<hD>^{-s} Op_h(a) <hD>^s f(x) = \int e^{i<\xi,x>} b(x,h\xi; h) \hat f(\xi) d\xi \]
provided we take
\[b(x,\xi;h) = \frac{1}{(2\pi h)^n}\int\int e^{\frac{i}{h}<x-y,\xi-\eta>}<\eta>^{-s}a(y,\xi)<\xi>^s dy d\eta\]
We now need to check that $b(x,\xi; h)$ does indeed satisfy the desired estimates. Make a change of variables and integrate by parts we get that
\[b(x,\xi;h) = \frac{1}{(2\pi)^n}\int\int e^{i<y,\eta>}(I - \Delta_\eta)^{n+1}\frac{<\xi-\eta>^{-s}<\xi>^s}{(1 + |\eta|^2)^{n+3}}\frac{(I -\Delta_y)^{n+3}a(x-hy,\xi;h)}{(1+|y|^2)^{n+1}}dyd\eta\]
By Peetre inequality, $<\xi-\eta>^{-s}<\xi>^s\leq <\eta>^s$. Since $|s|\leq 2$ this means that $\frac{<\xi-\eta>^{-s}<\xi>^s}{(1 + |\eta|^2)^{n+3}}\leq \frac{1}{(1 + |\eta|^2)^{n+1}}$ and is therefore integrable. It is easily seen by straight forward computation (or taking the logarithm then differentiate) that for all multi-indices $\beta$
\[\partial_\xi^\beta(I - \Delta_\eta)^{n+1}\frac{<\xi-\eta>^{-s}<\xi>^s}{(1 + |\eta|^2)^{n+3}} \leq C_{n,\beta} \frac{<\xi-\eta>^{-s}<\xi>^s}{(1 + |\eta|^2)^{n+3}}\]
Using this fact and Peetre inequality one sees that $b(x,\xi;h)$ satisfies for all multi-indices $\alpha,\beta$, 
\[p_{\alpha,\beta}^\sigma(b) \leq C_{\alpha,\beta}\sum\limits_{|\alpha'|,|\beta'|\leq n+3 + |\alpha|+|\beta|} p_{\alpha',\beta'}^\sigma(a)\] 

$\square$\\
We now derive some weighted space estimates for operators of class $Op_h S^0_\sigma$
\begin{lemma}
\label{weighted estimates}
Let $a\in S_\sigma^0$ with $0\leq \sigma\leq 1/2$. Then $Op_h(a)$ is bounded $L^2_\delta\to L^2_\delta$ for all $|\delta|\leq 2$. Furthermore, there exists a constant $C$ depending on dimension only such that
\[\|Op_h(a)\|_{L^2_{\delta}\to L^2_{\delta}}\leq C\sum\limits_{\alpha,\beta\leq k(n) + 2n+4} p_{\alpha,\beta}^\sigma(a)\]
for all $|\delta|\leq 2$, $0<h\leq 1$.
\end{lemma}
\noindent{\bf Proof}\\
If $\delta = -2$, define $Tf(x) = <x>^{-2} A(<x>^2 f)$. For $f \in {\cal S}$ one has
\begin{eqnarray*}
Tf(x) &=& (2\pi)^{-n}\int e^{ix\cdot\xi} a(x,h\xi) <x>^{-2}(I - \Delta_\xi)\hat f(\xi)d\xi\\
&=& (2\pi)^{-n}\int (I - \Delta_\xi)(e^{ix\cdot\xi} a(x,h\xi) <x>^{-2})\hat f(\xi)d\xi\\
\end{eqnarray*}
It is easily seen that $T$ is a semi-classical $\Psi DO$ with symbol $a'(x,\xi;h)\in S^0_\sigma$ and semi-norms of $a'(x,\xi;h)$ are bounded by
\[ p_{\alpha,\beta}^\sigma(a') \leq C_{\alpha,\beta}\sum_{|\alpha'|,|\beta'| \leq |\alpha|+|\beta| + 2}p_{\alpha',\beta'}^\sigma(a)\]
Use the above semi-norm estimate and apply proposition \ref{operator by symbol} to $a'$ we get the desired estimate for $\delta = -2$. To get the estimate for $Op_h(a)$ acting on $L^2_{2}$, we consider its adjoint $Op_h(a)^*$ acting on $L^2_{-2}$. Lemma \ref{adjoint symbol estimate} shows that $Op_h(a)^* = Op_h(a^*)$ with $a^*\in S^0_\sigma$ satisfying  
\[p_{\alpha,\beta}^\sigma (a^*) \leq C_{\alpha,\beta}\sum\limits_{|\alpha'|,|\beta'|\leq |\alpha| + |\beta| + 2n+2}p_{\alpha',\beta'}(a)\]
Apply the result we already have for $\delta = -2$ to $Op_h(a^*) = Op_h(a)^*$ acting on $L^2_{-2}$ in conjunction with the above estimate we get
\[\|Op_h(a)^*\|_{L^2_{-2}\to L^2_{-2}} \leq C\sum\limits_{\alpha,\beta\leq k(n) + 2n+4} p_{\alpha,\beta}^\sigma(a)\]
Now an argument using the duality between $L^2_2$ and $L^2_{-2}$ gives the estimate for $\delta =2$. Using interpolation we get the estimate for all $\delta\in[-2,2]$
$\square$\\
Using this lemma along with lemma \ref{conjugate} one can even obtain estimates for operators acting on weighted semi-classical sobolev spaces (see for example \cite{salo}).
\begin{lemma}
\label{weighted sobolev estimate}
Let $a\in S_\sigma^0$ with $0\leq \sigma\leq 1/2$. Then $Op_h(a)$ is bounded $H^s_{\delta,h}\to H^s_{\delta,h}$ for all $|\delta|\leq 2$. Furthermore, there exists a constant $C$ such that
\[\|Op_h(a)\|_{L^2_{\delta}\to L^2_{\delta}}\leq C\sum\limits_{\alpha,\beta\leq k(n) + 2n+4} p_{\alpha,\beta}^\sigma(a)\]
for all $|\delta|\leq 2$, $|s|\leq 2$ $0<h\leq 1$.
\end{lemma}
\noindent{\bf Proof}\\
By lemma \ref{conjugate}, the semi-norms for the symbol of $<hD>^s A <hD>^{-s}$ can be bounded by semi-norms of $a$. So apply lemma \ref{weighted estimates} to the $Op_h(S_\sigma^0)$ operator $<hD>^s A <hD>^{-s}$ we get
\begin{eqnarray}
\label{auxilary}
\|<x>^\delta <hD>^s Af\|_{L^2}\leq C \|<x>^\delta <hD>^s f\|_{L^2}\sum\limits_{|\alpha|,|\beta|\leq k(n) + 3n+7} p_{\alpha,\beta}^\sigma(a)
\end{eqnarray}
As shown in \cite{salo}, there exists a constant depending on dimension only such that the following inequality holds for all $f\in{\cal S}$ and $|\delta|\leq 2$
\[\|<hD>^2 <x>^\delta f\|_{L^2} \leq C\|<x>^\delta <hD>^2 f\|_{L^2}\]
Apply this inequality to $Af$ combined with (\ref{auxilary}) we get as in \cite{salo}
\begin{eqnarray*}
\|<hD>^2 <x>^\delta f\|_{L^2} &\leq& C(\sum\limits_{|\alpha|,|\beta|\leq k(n) + 3n+7} p_{\alpha,\beta}^\sigma(a))(\|<x>^\delta <hD>^2 f\|_{L^2})\\
&\leq&C(\sum\limits_{|\alpha|,|\beta|\leq k(n) + 3n+7} p_{\alpha,\beta}^\sigma(a))(\sum\limits_{|\alpha|\leq 2} \|<x>^\delta(hD)^\alpha f\|_{L^2})\\
&\leq&C(\sum\limits_{|\alpha|,|\beta|\leq k(n) + 3n+7} p_{\alpha,\beta}^\sigma(a))(\sum\limits_{|\alpha|\leq 2} \|(hD)^\alpha (<x>^\delta f)\|_{L^2})\\
&\leq&C(\sum\limits_{|\alpha|,|\beta|\leq k(n) + 3n+7} p_{\alpha,\beta}^\sigma(a))(\|<hD>^2<x>^\delta f\|_{L^2})
\end{eqnarray*}
So we have proven the inequality for the case $s = 0,2$. Moving to the Fourier side we see that this is equivalent to having for all $f\in{\cal S}$ satisfy
\[\|<h\xi>^2<D_\xi>^\delta \hat A <D_\xi>^{-\delta}f\|_{L^2} \leq C\{\sum\limits_{|\alpha|,|\beta|\leq k(n) + 3n+7} p_{\alpha,\beta}^\sigma(a)\}\|<h\xi>^2f\|_{L^2}\]
\[\|<D_\xi>^\delta \hat A <D_\xi>^{-\delta}f\|_{L^2} \leq C\{\sum\limits_{|\alpha|,|\beta|\leq k(n) + 3n+7} p_{\alpha,\beta}^\sigma(a)\}\|f\|_{L^2}\]
where $\hat A g := \widehat{A \check g}$. Interpolate the norm of the operator $<D_\xi>^\delta \hat A <D_\xi>^{-\delta}$ between the weighted spaces $<h\xi>^2$ and $<h\xi>^0$ gives the desired result for all $s\in[0,2]$. The case of $s\in [-2,0]$ can be done by similar duality argument used in the previous lemma.
$\square$\\

Suppose that $a\in S_\sigma^0$ satisfies $1/a\in S_\sigma^0$. It is well known that operators $Op_h(a)$ associated to such symbols are invertible provided $h > 0$ is taken to be smaller than some $h_0>0$ with $h_0$ depending on the chosen symbol. The next lemma addresses the question of when the $h_0$ can be taken uniformly for a given family of such symbols. More precisely
\begin{lemma}
\label{uniform invertibility}
For all $M>0$ and $\sigma < \frac{1}{2n+4}$ there exists an $h_0 >0$ such that for all symbols $a\in S_\sigma^0$ satisfying
\[\sum\limits_{|\alpha|,|\beta|\leq 4n+8 + k(n)} p^\sigma_{\alpha,\beta}(a) \leq M\ \ \ \sum\limits_{|\alpha|,|\beta|\leq 4n+8 + k(n)} p^\sigma_{\alpha,\beta}(\frac{1}{a}) \leq M\]

$Op_h(a)$ is invertible on $L^2_\delta$ for all $|\delta|\leq 2$ and $0<h \leq h_0$. Furthermore the norm of the inverse is uniformly bounded
\[\|Op_h(a)^{-1}\|_{L^2_\delta \to L^2_{\delta}}\leq C(M)\]
here $C(M)$ depends only on $M$ and the dimension. 
\end{lemma}
\noindent{\bf Proof}\\
By lemma \ref{multiplication remainder} $Op_h(a)Op_h(1/a) = I + h^{1-(2n + 4)\sigma}Op_h(m)$ with $m(x,\xi;h)\in S^0_\sigma$ satisfying the semi-norm estimates

\[p_{\alpha,\beta}^\sigma(m) \leq C_{\alpha,\sigma} (\sum\limits_{|\alpha'|,|\beta'|\leq |\alpha|+|\beta| + 2n+4}p_{\alpha',\beta'}^\sigma(a) + \sum\limits_{|\alpha'|,|\beta'|\leq |\alpha|+|\beta| + 2n+4}p_{\alpha',\beta'}^\sigma(a^{-1}))\]

Combining this and the weighted space estimates of lemma \ref{weighted estimates} we get
\[\|Op_h(m)\|_{L^2_\delta \to L^2_\delta} \leq C \sum\limits_{\alpha,\beta\leq 4n + 8 + k(n)}(p_{\alpha,\beta}^\sigma(a) +p_{\alpha,\beta}^\sigma(b))\leq CM\]
for all $|\delta|\leq 2$. Since $\sigma < \frac{1}{2n+4}$, we can pick $h_0>0$ such that $h_0^{1-(2n + 4)\sigma}CM \leq 1/2$. With this choice one sees that for all $0<h\leq h_0$, $(I + h^{1-(2n + 4)\sigma} Op_h(m))$ is invertible with 
\[\|(I + h^{1-(2n + 4)\sigma} Op_h(m))^{-1}\|_{L^2_\delta\to L^2_\delta} \leq 2\]
for all $|\delta|\leq 2$. So $Op_h(a)$ has a right inverse for $h \leq h_0$ that has norm bounded by $C(M)$. The exact same argument applied to $Op_h(1/a)Op_h(a)$ implies the existence of a left inverse of $Op_h(a)$ with norm bounded by $C(M)$. So there exists an $C >0$, $h_0 >0$ depending on dimension and $M$ only such that $Op_h(a)$ is bounded invertible from $L^2_\delta \to L^2_\delta$ with norm $\|Op_h(a)^{-1}\|_{L^2_\delta \to L^2_\delta} \leq CM$ for all $|\delta|\leq 2$ and $h \leq h_0$.
$\square$\\
\end{section}

\begin{section}{Properties of Complex Geometric Optic Solutions for Magnetic Schr\"odinger Equation}
Given an $M >0$, $R>0$ $p>n$ consider the family compactly supported vector fields 
\[{\cal W}(M,R) := \{W\in W^{1,p}(\R^n;\R^n) \mid \|W\|_{W^{1,p}}\leq M, \|div W\|_{L^\infty}\leq M\  and\  supp(W)\subset B_R\}\] 
and the family of compactly supported electric potentials 
\[{\cal Q}(M,R) := \{q\in L^\infty\mid \|q\|_{L^\infty}\leq M, supp(q)\subset B_R\}\] 
In this section we show that the CGO solutions to the equation $H_{W,q}u = 0$ has remainder decaying uniformly for all $W\in {\cal W}(M,R)$, $q\in {\cal Q}(M,R)$. More precisely,
\begin{proposition}
\label{existence of CGO}
Let $\sigma_0,\theta$ be positive numbers satisfying $\sigma_0 + \theta < \frac{1}{4n+6}$. For all $M>0$ and $R>0$ there exists constants $C>0$, $h_0 >0$, and $\epsilon >0$ depending on dimension, $p$, $\sigma_0$, $\theta$, $M$ and $R$ only such that for all $(W,q)\in {\cal W}(M,R)\times{\cal Q}(M,R)$, $\zeta \in \C^n$ with $\zeta\cdot\zeta = 0$ and $|\zeta|\geq \frac{1}{h_0}$, there exists solutions to $H_{W,q}u = 0$ in $\R^n$ of the form
\[u(x,\zeta) = e^{i\zeta\cdot x}(e^{i\chi_{|\zeta|}\phi^\sharp} + r(x,|\zeta|)),\ \ \ \ \ \|r(\cdot,|\zeta|)\|_{H^t_\delta}\leq C |\zeta|^{t-\epsilon},\ \ \ \ t\in[0,2]\]
where $\phi^\sharp$ is defined by
\[\phi^\sharp(x) = \frac{-\sqrt{2}\zeta}{2\pi |\zeta|}\cdot\int_{\R^2}\frac{ W^\sharp(x-y_1\frac{\sqrt{2}Re(\zeta)}{|\zeta|} -y_2\frac{\sqrt{2}Im(\zeta)}{|\zeta|})}{y_1+iy_2}\]
Here $\chi_{|\zeta|} = \chi(x/|\zeta|^{\theta})$ with $\chi(x)$ being a smooth function supported in the unit ball and is $1$ near zero, $W^\sharp$ is the convolution of $W$ with the mollifier defined by 
\[W^\sharp(x) = \int |\zeta|^{n\sigma_0}\chi(y|\zeta|^{\sigma_0}) W(x - y) dy\]
\end{proposition}
The next subsection proves some facts about the particular semi-classical symbols that we will be working with.

\begin{subsection}{Semi-Classical Symbols Arising from Magnetic Field}
In this section we define some symbols and prove some estimates for their semi-norms. For the motivation of these definitions see \cite{salo}. If $\zeta\in\C^n$ satisfies $\zeta\cdot\zeta = 0$, then we have $\zeta = \mu/h$ with $h = \sqrt{2}/|\zeta|$ and $\mu = \gamma_1 + i\gamma_2$ where $\gamma_1, \gamma_2\in\R^n$ are unit vectors satisfying $\gamma_1\perp\gamma_2$. Let $\sigma_0$ and $\sigma$ be positive number satisfying $0 <\sigma_0 <\sigma  < \frac{1}{4n + 6}$ and define $\theta = \sigma - \sigma_0$. For a given compactly supported magnetic potential $W\in {\cal W}(M,R)$ we decompose $W = W^\sharp + W^\flat$ where
\[W^\sharp(x) = \int \frac{1}{h^{n\sigma_0}}\chi(\frac{y}{h^{\sigma_0}}) W(x - y) dy\]
With these notations we define the nonsmooth symbol $r(x,\xi) = W(x)\cdot(\xi + \mu)$ and its smooth approximation $r^\sharp = W^\sharp(x)\cdot(\xi + \mu)$. Notice now that with this definition, 
\[|\partial_x^\alpha r^\sharp(x,\xi)| \leq C_{\alpha} \|W\|_{L^\infty}h^{-|\alpha|\sigma_0}(|\xi|+1)\leq C_{\alpha} Mh^{-|\alpha|\sigma_0}(|\xi|+1)\] Finally we define the elliptic symbol $q(\xi) = (\xi + \gamma_1)^2 -1 + 2i\gamma_2\cdot \xi$.\\
For $\epsilon >0$ we will consider the neighbourhood
\[U(\epsilon) = \{\xi\in\R^n\mid 1-\epsilon < |\xi + \gamma_1| < 1 + \epsilon, |\xi\cdot \gamma_2| <\epsilon\}\]
and introduce a smooth cutoff $\psi$ with $\psi =1$ in $U(\epsilon/4)$ and $\psi = 0$ outside of $U(\epsilon/2)$. Define the function
\begin{eqnarray}
\label{definition of w}
w(x,\xi) := \frac{-1}{2\pi}\int_{\R^n}\frac{1}{y_1 + iy_2}\psi(\xi) r^\sharp(x-y_1(\xi + \gamma_1) - y_2\gamma_2,\xi)dy_1dy_2
\end{eqnarray}
Then by lemma \ref{local transport estimate} it is a $C^\infty$ function that solves
\[(\xi + \mu)\cdot\nabla_x w = -\psi(\xi) r^\sharp(x,\xi)\]
and satisfies the estimates
\begin{eqnarray}
\label{estimates for w}
|\partial_x^\alpha \partial_\xi^\beta w(x,\xi; h)|\leq C_{\alpha,\beta,R}Mh^{-\sigma_0 |\alpha + \beta|}<x>^{|\beta|-1}
\end{eqnarray}
We see here that $w(x,\xi;h)$ is not quite a symbol of class $S^0_{\sigma_0}$ since it has growth in $x$ if we take more than 1 derivatives in $\xi$. To remedy this problem we will introduce another cutoff this time in the $x$ variable. Let $\chi$ be a compactly supported smooth function that is identically $1$ in $B_R$ and define
\[\phi(x,\xi) = \chi(h^\theta x) w(x,\xi)\]
where $\theta = \sigma - \sigma_0$. With $\theta$ chosen this way (\ref{estimates for w}) shows that $\phi$ satisfies
\begin{eqnarray}
\label{seminorms of phi}
|\partial_x^\alpha \partial_\xi^\beta \phi(x,\xi)| \leq C_{\alpha,\beta,R}M<x>^{-1}h^{-\sigma|\alpha+\beta|}
\end{eqnarray}
Therefore $\phi$ and $<x>\phi$ are in $S_\sigma^0$ with semi-norms $p_{\alpha,\beta}^\sigma(\phi)$ and $p_{\alpha,\beta}^\sigma(<x>\phi)$, bounded above by $C_{\alpha,\beta,R}M$. Furthermore if we define symbols $a$, $b$ by
\begin{eqnarray}
\label{definition of a, b}
a = e^{i\phi}\ \ \ b = a + 2h\frac{1-\psi(\xi)}{q(\xi)}e^{i\phi} r^\sharp
\end{eqnarray}
then $a$ and $b$ are in $S_\sigma^0$ with semi-norm estimate
\begin{eqnarray}
\label{estimates on a,b}
p_{\alpha,\beta}^\sigma(b) +p_{\alpha,\beta}^\sigma(a) \leq C_{\alpha,\beta,R}Me^M 
\end{eqnarray}
Note that although the notation does not explicitly state the dependence of $a$ and $b$ on the given magnetic field $W$, we must keep in mind that there is indeed a nontrivial dependence. The following lemma states a "uniform invertibility" result for $Op_h(a)$ and $Op_h(b)$.
\begin{lemma}
\label{invertibility of a,b}
There exists an $h_0>0$ such that for all $W\in {\cal W}(M,R)$, the symbol $a(x,\xi;h)$ arising from $W$ as in (\ref{definition of a, b}) has invertible semiclassical quantization $Op_h(a)$ for all $h\leq h_0$. More precisely, for all $h \leq h_0$, $Op_h(a)$ is invertible from $L^2_\delta \to L^2_\delta$ with norm
\[\|Op_h(a)^{-1}\|_{L^2_\delta\to L^2_\delta}\leq C_RMe^M\]
for all $|\delta|\leq 2$. The same holds for $b(x,\xi;h) = a(x,\xi;h) + 2h\frac{1-\psi(\xi)}{q(\xi)}e^{i\phi} r^\sharp $.
\end{lemma}
\noindent{\bf Proof}\\
Observe that $1/a = e^{-i\phi}$ so it satisfies the same semi-norm estimates as $a$. So by (\ref{estimates on a,b}) there exists an $M'$ such that all symbols $a$ arising from some $W\in {\cal W}(M,R)$ satisfies
\[\sum\limits_{|\alpha|,|\beta|\leq 4n+8 + k(n)} p^\sigma_{\alpha,\beta}(a) \leq M'\ \ \ \sum\limits_{|\alpha|,|\beta|\leq 4n+8 + k(n)} p^\sigma_{\alpha,\beta}(\frac{1}{a}) \leq M'\]
Lemma \ref{uniform invertibility} then gives an $h_0 >0$ such that all quantizations $Op_h(a): L^2_\delta \to L^2_\delta$ are invertible as long as $h\leq h_0$. Furthermore the norm of the inverse is bounded by
\[\|Op_h(a)^{-1}\|_{L^2_\delta\to L^2_\delta} \leq 2M'\]
for all $|\delta|\leq 2$. Moving on to $b$, observe that 
\[Op_h(b) = Op_h(a) + 2hOp_h(\frac{1-\psi(\xi)}{q(\xi)}e^{i\phi} r^\sharp)\] 
Obviously, $\frac{1-\psi(\xi)}{q(\xi)}e^{i\phi} r^\sharp\in S_\sigma^0$ with 
\[p_{\alpha,\beta}^\sigma(\frac{1-\psi(\xi)}{q(\xi)}e^{i\phi} r^\sharp)\leq C_{\alpha,\beta,R}Me^M\] 
for all $W\in {\cal W}(M,R)$. Therefore, by lemma \ref{weighted estimates} there exists a constant (WLOG) $M'$ such that 
\[\|Op_h(\frac{1-\psi(\xi)}{q(\xi)}e^{i\phi} r^\sharp\|_{L^2_\delta\to L^2_\delta}\leq M'\]
for all $h >0$, $|\delta|\leq 2$ and $W\in{\cal W}(M,R)$. Choose $h_0>0$ such that $2 h_0M'^2\leq \frac{1}{2}$ and a simple argument involving Neumann series shows that $Op_h(b)$ is invertible from $L^2_\delta\to L^2_\delta$ with norm
\[\|Op_h(b)^{-1}\|_{L^2_\delta\to L^2_\delta} \leq4M'\]
for all $h \leq h_0$ and $|\delta|\leq 2$.
$\square$\\
We conclude this subsection by introducing two other operators which will be key in the construction of CGO. Consider the symbol
\[h^{1+\theta}\{2e^{i\phi}[(\xi + \mu)\cdot w\nabla \chi(h^\theta x)]\}  + h^{2-2\sigma}\{h^{2\sigma}[\Delta_x a + 2W^\sharp\cdot D_x a]\}\] 
Due to estimate (\ref{seminorms of phi}) and the fact that the $w(x,\xi) = 0$ for $|\xi|\geq 1$, if we take $\epsilon = min\{\theta, 1-2\sigma\}$ this symbol can be written as $h^{1 + \epsilon} r_0$ with $r_0\in S^0_\sigma$ and $<x>r_0\in S_\sigma^0$. It can easily be checked from the definition that $p_{\alpha,\beta}^\sigma(<x>r_0) \leq C_{\alpha,\beta,R}Me^M$ for all $W\in {\cal W}(R,M)$. Therefore, by lemma \ref{weighted estimates} this implies that 
\begin{eqnarray}
\label{estimate on the operator R0}
\|Op_h(r_0)\|_{L^2_\delta\to L^2_{\delta+1}}\leq CMe^M\ \  for\ all \ \ \ |\delta| \leq 2
\end{eqnarray}
Now consider the operator $T$ and its inverse
\begin{eqnarray*}
\label{definition of T}
T = B\Delta_\zeta A^{-1}\Delta_\zeta^{-1},\ \ \ T^{-1} = \Delta_\zeta A\Delta_\zeta^{-1}B^{-1}
\end{eqnarray*}

By \cite{salo} they can be written as
\[T = I + 2W^\sharp \cdot D_\zeta\Delta_\zeta^{-1} - h^{-1 + \epsilon}Op_h(r_0) A^{-1}\Delta_\zeta^{-1}\]
\[T^{-1} = I - 2W^\sharp\cdot D_\zeta A\Delta_\zeta^{-1}B^{-1} + h^{-1+\epsilon}Op_h(r_0)\Delta_\zeta^{-1}B^{-1}\]
It is now easy to see by analyzing the operators term by term that they are both bounded operators from $L^2_{\delta+1}\to L^2_{\delta+1}$ with norm depending only on $M$. For example, lets take 
\[2W^\sharp\cdot D_\zeta A\Delta_\zeta^{-1}B^{-1}: L^2_{\delta+1}\to L^2_{\delta +1}\] 
By lemma \ref{invertibility of a,b} the operator $B^{-1}: L^2_{\delta+1}\to L^2_{\delta +1}$ with norm $CMe^M$. Classical results for the $\Delta_\zeta^{-1}$ operator shows that $\Delta_\zeta^{-1} : L^2_{\delta + 1}\to H^1_{\delta,|\zeta|^{-1}}$ with norm less than $C|\zeta|^{-1}$ with $C$ depending only on $\delta$ and the dimension. Lemma \ref{weighted sobolev estimate} combined with inequality (\ref{estimates on a,b}) shows that $A : H^1_{\delta,|\zeta|^{-1}}\to H^1_{\delta,|\zeta|^{-1}}$ with norm $CMe^M$. Since $W^\sharp$ is compactly supported with support independent of $\zeta$, $W^\sharp\cdot D_\zeta :  H^1_{\delta,|\zeta|^{-1}}\to L^2_{\delta +1}$ with norm less than $M|\zeta|$. So we see that $2W^\sharp\cdot D_\zeta A\Delta_\zeta^{-1}B^{-1}: L^2_{\delta+1}\to L^2_{\delta +1}$ with norm at most $CMe^M$ with $C$ depending on dimension and $\delta$ only. The rest of the terms can be handled in a similar way to give that
\begin{eqnarray}
\label{norm estimate for T}
\|T\|_{L^2_{\delta+1}\to L^2_{\delta+1}}\leq CMe^M\ \ \ \|T^{-1}\|_{L^2_{\delta+1}\to L^2_{\delta+1}}\leq CMe^M
\end{eqnarray}
for $-1<\delta<0$.

\end{subsection}
\begin{subsection}{Invertibility of the $(\Delta_\zeta + 2W\cdot D_\zeta + q)$ Operator}
In this section we prove the following theorem which plays a critical role in the construction of CGO solutions
\begin{proposition}
\label{solution to the remainder equation}
Fix $\delta\in (-1,0)$. For all $M>0$, $R>0$ there exists an $C>0$ and $h_0>0$ such that for all $W\in {\cal W}(M,R)$, $q\in {\cal Q}(M,R)$, and $f\in L^{2}_{\delta+1}$ the equation
\[(\Delta_\zeta + 2W\cdot D_\zeta + q)u = f\]
has a unique solution $u\in H^1_{\delta}$ for all $\zeta\in\C^n$ satisfying $\zeta\cdot\zeta = 0$ and $|\zeta|\geq \frac{1}{h_0}$. Furthermore the following estimates holds for $u$ with constants depending only on $M$ and $R$ but not on the choice of $(W,q)\in {\cal W}(M,R)\times{\cal Q}(M,R)$:
\[\|u\|_{H^t_\delta}\leq C\|f\|_{L^2_{\delta + 1}}|\zeta|^{t-1}\ \ \ \ t\in[0,2]\]
\end{proposition}
\noindent{\bf Proof}\\
As in \cite{salo} we seek solutions in the form of $u = \Delta_\zeta^{-1}v$. Following \cite{salo} we see that $v$ satisfies
\[(I + h^{-1+\epsilon}Op_h(r_0) A^{-1}\Delta_\zeta^{-1}T^{-1} + 2 W^{\flat}\cdot D_\zeta \Delta_\zeta^{-1}T^{-1} + q\Delta_\zeta^{-1}T^{-1})Tv = f\]
First we will show the invertibility of the operator 
\[(I + h^{-1+\epsilon}Op_h(r_0) A^{-1}\Delta_\zeta^{-1}T^{-1} + 2 W^{\flat}\cdot D_\zeta \Delta_\zeta^{-1}T^{-1} + q\Delta_\zeta^{-1}T^{-1})\]
from $L^2_{\delta + 1}\to L^2_{\delta + 1}$. By (\ref{estimate on the operator R0}), (\ref{norm estimate for T}), and lemma \ref{invertibility of a,b}, there exists an $h_0>0$ such that for all $h\leq h_0$,  
\[\|h^{-1+\epsilon}Op_h(r_0) A^{-1}\Delta_\zeta^{-1}T^{-1}\|_{L^2_{\delta + 1}\to L^2_{\delta + 1}}\leq Ch^\epsilon\]
with $C$ depending only on $M$.
To obtain the same result for $W^{\flat}\cdot D_\zeta \Delta_\zeta^{-1}T^{-1}$, we observe that $W\in W^{1,p}_c(B_R;\R^n)\hookrightarrow C_c^{0,1-n/p}(B_R;\R^n)$ with norm $\|W\|_{C^{0,1-n/p}(B_R;\R^n)}\leq C\|W\|_{W^{1,p}(B_R;\R^n)}$ so that $\|W^\flat\|_{L^\infty}\leq CMh^{1-n/p}$. Therefore 
\[\|W^{\flat}\cdot D_\zeta \Delta_\zeta^{-1}T^{-1}\|_{L^2_{\delta + 1}\to L^2_{\delta + 1}}\leq C h^{1-n/p}\]
The next term can be handled easily by using the same analysis and the fact that $supp(q)\subset B_R$. So by a Neumann series argument we see that there exists an $h_1 >0$ such that for all $h\leq h_1$ the operator
\[(I + h^{-1+\epsilon}Op_h(r_0) A^{-1}\Delta_\zeta^{-1}T^{-1} + 2 W^{\flat}\cdot D_\zeta \Delta_\zeta^{-1}T^{-1} + q\Delta_\zeta^{-1}T^{-1})\]
is invertible from $L^2_{\delta + 1}\to L^2_{\delta + 1}$ with norm of inverse less than $2$. Due to (\ref{norm estimate for T}), $T$ is invertible with norm of the inverse bounded by a constant depending only on $M$. Therefore we can write
\[v = T^{-1}(I + h^{-1+\epsilon}Op_h(r_0) A^{-1}\Delta_\zeta^{-1}T^{-1} + 2 W^{\flat}\cdot D_\zeta \Delta_\zeta^{-1}T^{-1} + q\Delta_\zeta^{-1}T^{-1})^{-1}f\]
with the estimate $\|v\|_{L^2_{\delta +1}}\leq C\|f\|_{L^2_{\delta+1}}$. The result now follows by writing $u = \Delta_\zeta^{-1}v$ and use the estimate we have for the operator $\Delta_\zeta^{-1}$.$\square$
\end{subsection}
\begin{subsection}{Proof of Proposition \ref{existence of CGO}}
We prove theorem \ref{existence of CGO} in this subsection. The proof is identical to that which appeared in \cite{salo} and \cite{sun uhlmann} except that we have to check that the decay of the remainder term is uniform for all $(W,q)\in {\cal W}(M,R)\times {\cal Q}(M,R)$.
Plug $u = e^{i\zeta\cdot x}(e^{i\phi^\sharp}+ r)$ into $H_{W,q}u$ we see that $r(x,|\zeta|)$ satisfies
\begin{eqnarray}
\label{equation for which remainder satisfies}
(\Delta_\zeta + 2W\cdot D_\zeta +G)r = -f
\end{eqnarray}
where $G = W^2 + D\cdot W + q\in L^\infty$ and 
\begin{eqnarray*}
f = (\Delta_\zeta + 2W\cdot D_\zeta +G)e^{i\chi_{|\zeta|}\phi^\sharp}= e^{i\chi_{|\zeta|}\phi^\sharp}\big\lbrack i\chi_{|\zeta| } \Delta \phi^\sharp + 2iD\chi_{|\zeta| }\cdot D\phi^\sharp + i\phi^\sharp\Delta\chi_{|\zeta| }\\
(\chi_{|\zeta| }\nabla \phi^\sharp +\phi^\sharp \nabla\chi_{|\zeta| })^2 + 2\zeta\cdot(\nabla\chi_{|\zeta| })\phi^\sharp + 2\zeta\cdot(\nabla\phi^\sharp)\chi_{|\zeta| }\\
+ 2W\cdot(\nabla\chi_{|\zeta| })\phi^\sharp + 2W\cdot(\nabla\phi^\sharp)\chi_{|\zeta| } + 2W^\sharp\cdot\zeta + 2W^\flat\cdot\zeta +G \big\rbrack
\end{eqnarray*}
By the choice of $\phi^\sharp$ it solves the transport equation
\[2\zeta\cdot\nabla\phi^\sharp + 2W^\sharp\cdot\zeta = 0\]
Since for $\zeta$ large $W^\sharp = \chi_\zeta W^\sharp$ the above relation removes the two terms that is of order $|\zeta|$ from $f$. So in terms of $L^2_{\delta + 1}$ norms this is 
\begin{eqnarray}
\label{weighted norm of f}
\|f\|&\leq& C\big\lbrack ||i\chi_{|\zeta| } \Delta \phi^\sharp|| + ||\nabla\chi_{|\zeta| }\cdot \nabla\phi^\sharp|| + ||\phi^\sharp\Delta\chi_{|\zeta| }||
+|||\chi_{|\zeta| }\nabla \phi^\sharp|^2|| +|||\phi^\sharp \nabla\chi_{|\zeta| }|^2||\\\nonumber &+& |\zeta|^{1-\theta}||\phi^\sharp\nabla\chi(x|\zeta|^{-\theta})|| + ||W\cdot(\nabla\chi_{|\zeta| })\phi^\sharp|| +||W\cdot(\nabla\phi^\sharp)\chi_{|\zeta| }|| \\\nonumber &+& |\zeta|||W^\flat|| +||G|| \big\rbrack
\end{eqnarray}
Every term except for $|\zeta|||W^\flat||$ in the above expression are of order $C|\zeta|^{1-\epsilon}$ for some $\epsilon >0$ by our choice of $\delta, \sigma, \theta$ and lemma \ref{L-infinity bound} (see \cite{salo} for details). Assuming without loss of generality that $\epsilon \leq 1-n/p$, Sobolev embedding gives $W\in C^{0,\epsilon}_c(B_R)$ with $\|W\|_{C^{0,\epsilon}(B_R)} \leq M$. So by property of the mollifier we have that $\|W^\flat\|_{L^2_{\delta+1}}\leq C|\zeta|^{-\epsilon\theta }$ with $C$ depending only on $M$ and $R$. Therefore we conclude that
\[\|f\|_{L^2_{\delta+1}}\leq C|\zeta|^{1-\epsilon}\]
for some $\epsilon >0$. Apply lemma \ref{solution to the remainder equation} to the operator $(\Delta_\zeta + 2W\cdot D_\zeta +G)$ gives a $C>0$ and $h_0>0$ such that for all $|\zeta|\geq 1/h_0$ and $\zeta\cdot\zeta =0$ equation (\ref{equation for which remainder satisfies}) has a unique solution satisfying
\[\ \ \ \|r\|_{H^t_\delta}\leq C|\zeta|^{t-\epsilon},\ \ \ \ t\in[0,2]\]
and the proof is complete.$\square$
\end{subsection}

\end{section}
\end{section}
\begin{section}*{PART II - Full Data Estimate}

\begin{section}{Estimate for the Magnetic Field}

We wish to derive log-type estimates for the curl of the difference of the magnetic potential $W_1 - W_2$. The main result of this section is
\begin{theorem}
\label{curl estimate}
Let $\Omega$ be a bounded open subset of $\R^n$ with smooth boundary. For all $M>0$, there exists a $C> 0$, $\epsilon>0$ such that the estimate
\[\|I_{\Omega}d(W_1 - W_2)\|_{H^{-1}(\R^n)}\leq C\lbrace \|\Lambda_{W_1,q_1}-\Lambda_{W_2,q_2}\|_{\frac{1}{2},\frac{-1}{2}}^{1/2} + |log \|\Lambda_{W_1,q_1}-\Lambda_{W_2,q_2}\|_{\frac{1}{2},\frac{-1}{2}}|^{-\epsilon}\rbrace\]

holds for all $W_1, W_2\in W^{1,\infty}(\Omega)$ and $q_1,q_2\in L^\infty(\Omega)$ satisfying $\|W_l\|_{ W^{1,\infty}(\Omega)}\leq M$, $\|q_l\|_{L^\infty(\Omega)}\leq M$ ($l = 1,2$) and $W_1 = W_2$ on $\partial\Omega$. Here $I_\Omega(x)$ is the indicator function of $\Omega$.
\end{theorem}
Before we proceed with the proof of theorem \ref{curl estimate}, we will see that the $W_1\mid_{\partial\Omega} = W_2\mid_{\partial\Omega}$ condition can be slightly relaxed provided we assume more regularity about the magnetic potentials. The only reason we need to make the additional assumption is to ensure that the gauge transformation one needs to do to reduce to the case of theorem \ref{curl estimate} has sufficient regularity. Sharper results with less a-priori assumptions are possible but the technical details would obscure the main points of this exposition.
\begin{corollary}
\label{reduction}
Let $\Omega$ be a bounded open subset of $\R^n$ with smooth boundary. For all $M>0$, there exists a $C> 0$, $\epsilon>0$ such that the estimate
\[\|I_{\Omega}d(W_1 - W_2)\|_{H^{-1}(\R^n)}\leq C\lbrace \|\Lambda_{W_1,q_1}-\Lambda_{W_2,q_2}\|_{\frac{1}{2},\frac{-1}{2}}^{1/2} + |log \|\Lambda_{W_1,q_1}-\Lambda_{W_2,q_2}\|_{\frac{1}{2},\frac{-1}{2}}|^{-\epsilon}\rbrace\]

holds for all $W_1, W_2\in W^{2,\infty}(\Omega)$ and $q_1,q_2\in L^\infty(\Omega)$ satisfying $\|W_l\|_{ W^{2,\infty}(\Omega)}\leq M$, $\|q_l\|_{L^\infty(\Omega)}\leq M$ ($l = 1,2$) and $W_1\mid_D = W_2\mid_D$ on $\partial\Omega$. Here $I_\Omega(x)$ is the indicator function of $\Omega$.
\end{corollary}
\noindent{\bf Proof}\\
Suppose the statement is proven in the case when $W_1 = W_2$ on $\partial\Omega$. Then for arbitrary $W_1, W_2$ such that $W_1\mid_D = W_2\mid_D$ and $\|W_l\|_{W^{2,\infty}(\Omega)}\leq M$, one can construct potentials $p_1,p_2\in W^{2,\infty}(\Omega)\cap H^1_0(\Omega)$ satisfying 
\[\| p_l\|_{W^{2,\infty}(\Omega)}\leq CM,\ \ \ \frac{\partial p_l}{\partial\nu} = -W_l\mid_N\]
where $C$ depends on $\Omega$ only. Then $\|W_l + \nabla p_l\|_{W^{1,\infty}(\Omega)}\leq CM$ with
\[W_1 + \nabla p_1 = W_2 + \nabla p_2\ \ \ on\ \ \ \partial\Omega\]

Apply theorem \ref{curl estimate} and use the gauge invariance property of both the DN map and the $d$ operator we get
\begin{eqnarray*}
\|I_\Omega d(W_1 - W_2)\|_{H^{-1}(\R^n)} &=& \|I_\Omega d((W_1-\nabla p_1) - (W_2-\nabla p_2))\|_{H^{-1}(\R^n)}\\ &\leq&  C\lbrace \|\Lambda_{W_1-\nabla p_1,q_1}-\Lambda_{W_2 - \nabla p_2,q_2}\|_{\frac{1}{2},\frac{-1}{2}}^{1/2} + |log \|\Lambda_{W_1-\nabla p_1,q_1}-\Lambda_{W_2-\nabla p_2,q_2}\|_{\frac{1}{2},\frac{-1}{2}}|^{-\epsilon}\rbrace\\
&\leq& C\lbrace \|\Lambda_{W_1,q_1}-\Lambda_{W_2,q_2}\|_{\frac{1}{2},\frac{-1}{2}}^{1/2} + |log \|\Lambda_{W_1 ,q_1}-\Lambda_{W_2,q_2}\|_{\frac{1}{2},\frac{-1}{2}}|^{-\epsilon}\rbrace
\end{eqnarray*}
$\square$\\\\
\begin{subsection}{Extending Vector Fields to a Larger Domain}
First we need to prove a technical lemma about extending a $W^{1,\infty}(\Omega)$ vector field to a slightly larger domain\\
\begin{lemma}
\label{gluing functions}
Let $\Omega\subset\subset\tilde\Omega$. Suppose $u\in W^{1,p}(\Omega)$, $u'\in W^{1,p}(\tilde\Omega\backslash\Omega)$ and $u = u'$ on $\partial\Omega$. Then the function defined by
\begin{eqnarray*}
U(x) := {u(x)\ \ x\in\Omega \atopwithdelims \{  \} u'(x)\ \ x\in\tilde\Omega\backslash\Omega}\qquad
\end{eqnarray*}
is in $W^{1,p}(\tilde\Omega)$ and
\begin{eqnarray*}
\partial_j U(x) = {\partial_j u(x)\ \ x\in\Omega \atopwithdelims \{  \} \partial_j u'(x)\ \ x\in\tilde\Omega\backslash\Omega}\qquad
\end{eqnarray*}
Furthermore, $\|U\|_{W^{1,p}(\tilde\Omega)} = \|u\|_{W^{1,p}(\Omega)} +\|u'\|_{W^{1,p}(\tilde\Omega\backslash\Omega)}$
\end{lemma}
\noindent{\bf Proof}
Integrate by parts and check that the definition of weak derivative holds.$\square$
\begin{lemma}
\label{divergence extention}
Let $\Omega$ be a bounded domain in $\R^n$ that is compactly contained in $\tilde\Omega$. Suppose $W_1, W_2 \in W^{1,\infty}(\Omega)$  such that for $l = 1,2$ $\|W_l\|_{ W^{1,\infty}(\Omega)}\leq M$  and $W_1 = W_2$ on $\partial\Omega$. Then there exists extentions $\tilde W_1, \tilde W_2\in W^{1,\infty}_c(\tilde\Omega)$ such that
\[\|\tilde W_l\|_{W^{1,\infty}(\tilde\Omega)}\leq C(\tilde\Omega,\Omega)M\]
and $W_1 = W_2$ on $\tilde\Omega\backslash\Omega$. Here the constant depends only on $\Omega$ and $\tilde\Omega$ but not on $M$. 
\end{lemma}
\noindent{\bf Proof}\\
By standard extention theorem, there exists a $\tilde W_1\in W^{1,\infty}_c(\tilde\Omega)$ such that $\|W_1\|_{W^{1,\infty}_c(\tilde\Omega)} \leq C\|W_1\|_{W^{1,\infty}(\Omega)}$. For $\tilde W_2$ we define 
\[\tilde W_2:= {W_2\ \ x\in\Omega \atopwithdelims \{  \} \tilde W_1\ \ x\in\tilde\Omega\backslash\Omega}\qquad\]

Since $W_1 = W_2$ on $\partial\Omega$, by lemma \ref{gluing functions} we have $\tilde W_2\in W^{1,\infty}_c(\tilde\Omega)$ with $\|\tilde W_2\|_{W^{1,\infty}(\tilde\Omega)}\leq CM$ and 
\[\partial_j\tilde W_2 = {\partial_j W_2\ \ x\in\Omega \atopwithdelims \{  \} \partial_j\tilde W_1\ \ x\in\tilde\Omega\backslash\Omega}\qquad\]
Therefore we have that $\tilde W_2 \in W^{1,\infty}_c(\tilde\Omega)$ and $\|W_2\|_{W^{1,\infty}(\tilde\Omega)}\leq CM$$\square$
\end{subsection}

\begin{subsection}{Proof of Theorem \ref{curl estimate}}
Fix $j,k \in \{1...n\}$. For all $\xi$ such that $\xi_j\not= 0$ or $\xi_k\not= 0$, define $\tilde \gamma := \frac{e_j\xi_k - e_k\xi_j}{|e_j\xi_k - e_k\xi_j|}$.\\
Choose $\gamma\perp\tilde\gamma$, $\tilde\gamma\perp\xi$, $ |\gamma| = 1$  and define $\mu = \gamma+i\tilde\gamma$. For each $s\geq|\xi|$ let
\[\zeta_1(s) = \zeta_1 = -is\tilde\gamma + g(\xi,s)\gamma + \xi\]
\[\zeta_2(s) = \zeta_2 = is\tilde\gamma + g(\xi,s)\gamma -\xi\]
where $g(\xi,s) = \sqrt{s^2 -|\xi|^2}$. Note that $\zeta_l(s)$ can be written as $\zeta_l(s) = (-1)^{l}is\tilde \gamma + s\gamma_l'$ $(l = 1,2)$ provided we take 
\[\gamma'_1 = \frac{g(\xi,s)\gamma + \xi}{s},\ \ \ \gamma'_2 = \frac{g(\xi,s)\gamma - \xi}{s}\]
Let $\sigma_0>0$ and $\theta>0$ be chosen so that $0<\sigma_0 <\sigma_0+\theta<\frac{1}{4n+6}$. Given any $W_1, W_2 \in W^{1,\infty}(\Omega)$ such that $W_1 = W_2$ on $\partial\Omega$ and $\|W_l\|_{W^{1,\infty}(\Omega)}\leq M$,  extend them to compactly supported vector fields $\tilde W_1, \tilde W_2\in W^{1,p}_c(\tilde\Omega;\R^n)$ as described in lemma \ref{divergence extention} of the previous section. Therefore, there exists an $M'>0$ and $R>0$ such that $\tilde W_l \in {\cal W}(M',R)$ whenever $\|W_l\|_{W^{1,\infty}(\Omega)}\leq M$. Proposition \ref{existence of CGO} then gives a $C>0$, $h_0 >0$, and $ 0<\epsilon< \sigma_0(1-\frac{n}{p})$ such that for all $|\zeta_l|\geq\frac{1}{h_0}$ there exists solutions $u_l(x,\zeta_l)$ to $H_{\tilde W_l,q_l}u_l = 0$ in $\R^n$ that are of the form:
\[u_l(x,\zeta_l) = e^{i\zeta_l\cdot x}(e^{i\chi_{|\zeta_l|}\phi^{\sharp}_l} + r_l(x,\zeta_l))\ \ \ \|r_l(\cdot,\zeta_l)\|_{H^t_\delta}\leq C|\zeta_l|^{t-\epsilon}\ 
\ \ \ t\in[0,2]\ \ \ -1<\delta<0\]
\[\phi^\sharp_l(x) = -\frac{\gamma_l' +(-1)^li\tilde\gamma}{2\pi}\cdot\int_{\R^2}\frac{\tilde W_l^\sharp(x -\gamma_l'y_1 - (-1)^l\tilde\gamma y_2)}{y_1 + iy_2}\]
where $\chi_{|\zeta_l|} = \chi(x/|\zeta_l|^{\theta})$ and $\tilde W_l^\sharp = \int |\zeta_l|^{n\sigma_0}\chi(y|\zeta_l|^{\sigma_0}) \tilde W_l(x - y) dy$


Note here that $C$, $h_0$, and $\epsilon$ depend on dimension, $\Omega$, $\theta$, $\sigma_0$ and $M$ only but not on the choice of $(W_l,q_l)$ as long as $\|W_l\|_{W^{1,\infty}(\Omega)}\leq M$ and $\|q_l\|_{L^\infty}\leq M$. All constants in this section will have only the aforementioned dependence. In fact, the only reason we need the a-priori bound to be on the $W^{1,\infty}(\Omega)$ norm is so that the extention $\tilde W_l$ would have $L^\infty(\R^n)$ divergence which is required for the construction of CGO by proposition \ref{existence of CGO}. After the solutions have been constructed we only need the a-priori estimate to be in the $W^{1,p}(\Omega)$ norm for $p >n$. To demonstrate this fact we will only work with this norm in the proof below.
\begin{lemma}
\label{leading term estimate}
Let $u_l$ be the solutions constructed above. For all $s\geq 1/h_0$ such that $\chi_s(x) =1$ for all $x\in \tilde\Omega$ we have
\[\int_\Omega(W_1 - W_2)\cdot(u_1\nabla\bar u_2-\bar u_2\nabla u_1)dx = -i\int_\Omega e^{i(\phi^{\sharp}_1 - \bar\phi^{\sharp}_2)}e^{2i\xi\cdot x}(W_1 - W_2)\cdot(\bar\zeta_2+\zeta_1) + f_1(x,s) + f_2(x,s)dx\]
with $\|f_1(x,s)\|_{L^1(\Omega)} + \|f_2(x,s)\|_{L^1(\Omega)}\leq Cs^{1-\epsilon}$
\end{lemma}
\noindent{\bf Proof}\\
For all $x\in\tilde\Omega$ direct calculation gives 
\begin{eqnarray}
\label{leading term}
\nonumber \bar u_2 \nabla u_1 &=& \zeta_1 e^{2i\xi\cdot x}e^{i(\phi_1^\sharp -\bar\phi_2^\sharp)} + \zeta_1 e^{2i\xi\cdot x}\{r_1(x,s)e^{i\bar\phi_2^\sharp} + \bar r_2(x,s)e^{i\bar\phi_1^\sharp} + r_1 \bar r_2\}\\ &+& e^{2i\xi\cdot x}\{ e^{i(\phi_1^\sharp -\bar\phi_2^\sharp)}\nabla \phi_1^\sharp + ie^{i\phi_1^\sharp}\bar r_2\nabla\phi_1^\sharp + e^{-i\bar\phi_2^\sharp}\nabla r_1 + \bar r_2 \nabla r_1\}\\\nonumber
&=&\zeta_1 e^{2i\xi\cdot x}e^{i(\phi_1^\sharp -\bar\phi_2^\sharp)} + T_1 +T_2
\end{eqnarray}
where
\[T_1 := \zeta_1 e^{2i\xi\cdot x}\{r_1(x,s)e^{i\bar\phi_2^\sharp} + \bar r_2(x,s)e^{i\bar\phi_1^\sharp} + r_1 \bar r_2\}\]
\[T_2 := e^{2i\xi\cdot x}\{ e^{i(\phi_1^\sharp -\bar\phi_2^\sharp)}\nabla \phi_1^\sharp + ie^{i\phi_1^\sharp}\bar r_2\nabla\phi_1^\sharp + e^{-i\bar\phi_2^\sharp}\nabla r_1 + \bar r_2 \nabla r_1\}\]
Use lemma \ref{L-infinity bound} and the fact that $\tilde W_l^\sharp$ ($l = 1,2$) is compactly supported in $\tilde\Omega$ 
\[\|\phi_l^{\sharp}\|_{L^{\infty}(\Omega)}\leq C\|\tilde W^\sharp_l\|_{L^{\infty}(\tilde\Omega)}\leq CM\]
Use the fact that $\|r_l\|_{L^2_\delta} \leq Cs^{-\epsilon}$, one easily shows that $\|T_1\|_{L^1(\Omega)}\leq Cs^{1-\epsilon}$. The next term $T_2$ can be handled in a similar fashion. The only problematic part of $T_2$ are the terms involving $\nabla \phi_1^\sharp$. But by lemma \ref{L-infinity bound} we have
\[\|\phi_1^\sharp\|_{W^{1,\infty}(\Omega)}\leq C\|\tilde W_1^\sharp\|_{W^{1,\infty}(\tilde\Omega)} \leq Cs^{\sigma_0}\|\tilde W_1^\sharp\|_{L^\infty(\tilde\Omega)}\leq Cs^{\sigma_0}\|\tilde W_1^\sharp\|_{W^{1,p}(\tilde\Omega)} \leq CM's^{\sigma_0}\leq CM's^{1-\epsilon}\]
the last inequality comes from the fact that $\sigma_0 < \frac{1}{4n + 6}$ and $\epsilon < \sigma_0(1-\frac{n}{p})$. So arguing term by term in $T_2$ we get that $\|T_2\|_{L^1(\Omega)}\leq Cs^{1-\epsilon}$. Combining these observations into (\ref{leading term}) we get that 
\[\int_\Omega (W_1-W_2)\cdot\bar u_2 \nabla u_1 = \int_\Omega(W_1 - W_2)\cdot\zeta_1 e^{2i\xi\cdot x}e^{i(\phi_1^\sharp -\bar\phi_2^\sharp)} + \int_\Omega T_1 +\int_\Omega T_2\]
with the $L^1(\Omega)$ norm of $T_1, T_2$ controlled by $Cs^{1-\epsilon}$. Of course, similar calculation holds for the $u_1 \nabla \bar u_2$ and the proof is complete. 
$\square$


\begin{lemma}
\label{curl identity}
\[\int_\Omega(W_1 - W_2)\cdot(\bar\zeta_2 + \zeta_1)e^{i(\phi_1^\sharp - \bar\phi_2^\sharp)}e^{2i\xi\cdot x}dx = -\int_{\Omega}(W_1 - W_2)\cdot 2s\bar\mu e^{2i\xi\cdot x}dx+ G(\xi,s)\]
with $|G(\xi,s)|\leq C|s - g(\xi,s)| + s^{1-\epsilon} + s|\tilde\gamma - \frac{g(\xi,s)\tilde\gamma + \xi}{s}|^{1-\frac{n}{p}}$
\end{lemma}
\noindent{\bf Proof}\\
By the definition of $\zeta_1,\zeta_2$ we can write 
\begin{eqnarray}
\label{plus then minus}
\nonumber\int_\Omega(W_1 - W_2)\cdot(\bar\zeta_2 + \zeta_1)e^{i(\phi_1^\sharp - \bar\phi_2^\sharp)}e^{2i\xi\cdot x}dx &=& -\int_\Omega(W_1 - W_2)\cdot 2s\bar\mu e^{i(\psi_1 - \bar\psi_2)}e^{2i\xi\cdot x}dx\\\nonumber 
&+&  \int_\Omega(W_1 - W_2)\cdot 2s\bar\mu (e^{i(\psi_1 - \bar\psi_2)}-e^{i(\psi_1^\sharp - \bar\psi_2^\sharp)})e^{2i\xi\cdot x}dx\\\nonumber
&+&\int_\Omega(W_1 - W_2)\cdot 2s\bar\mu (e^{i(\psi_1^\sharp - \bar\psi_2^\sharp)}- e^{i(\phi_1^\sharp -\bar\phi_2^\sharp)})e^{2i\xi\cdot x}dx\\\nonumber
&+&\int_\Omega(W_1 - W_2)\cdot\gamma 2(s - g(\xi,s))e^{i(\phi_1^\sharp- \bar\phi_2^\sharp)}e^{2i\xi\cdot x}\\
&& := T_1 + T_2+T_3+T_4
\end{eqnarray}
where $\psi_1 := N^{-1}_{\bar\mu}(-\bar\mu\cdot \tilde W_1)$, $\psi_2 := N^{-1}_{\mu}(-\mu\cdot \tilde W_2)$ and $\psi_1^\sharp, \psi_2^\sharp$ are defined the same way with $\tilde W_l^\sharp$ replacing $\tilde W_l$. A simple calculation shows that $\psi_1(x) -\bar\psi_2(x) = N_{\bar\mu}^{-1}(-\bar\mu\cdot(\tilde W_1 - \tilde W_2))$. Recall that $\tilde W_1 = \tilde W_2$ on $\R^n\backslash\Omega$, so by lemma \ref{nonlinear fourier}
\[T_1 = -\int_{\R^n}(\tilde W_1 - \tilde W_2)\cdot 2s\bar\mu e^{i(\psi_1 - \bar\psi_2)}e^{2i\xi\cdot x}dx =  -\int_{\R^n}(\tilde W_1 - \tilde W_2)\cdot 2s\bar\mu e^{2i\xi\cdot x} = -\int_{\Omega}( W_1 - W_2)\cdot 2s\bar\mu e^{2i\xi\cdot x}\]
Now it remains to estimate the three other terms. Applying lemma \ref{approximation by mollifier}, \ref{exponential estimate}, and \ref{L-infinity bound} to $T_2$, $T_3$ and $T_4$ respectively, and use the fact that $W_c^{1,p}(\tilde\Omega)\hookrightarrow C^{0,1 - \frac{n}{p}}(\tilde\Omega)$ for $p>n$ we get that
\[|T_2| \leq Cs(\|\tilde W_1\|_{C^{0,1 - \frac{n}{p}}(\tilde\Omega)} + \|\tilde W_2\|_{C^{0,1 - \frac{n}{p}}(\tilde\Omega)})s^{-\sigma_0 (1-\frac{n}{p})}\leq Cs^{1-\epsilon}\]
\[|T_3|\leq Cs (\|\tilde W_1\|_{C^{0,1 - \frac{n}{p}}(\tilde\Omega)} + \|\tilde W_2\|_{C^{0,1 - \frac{n}{p}}(\tilde\Omega)}) |\tilde\gamma - \frac{g(\xi,s)\tilde\gamma + \xi}{s}|^{1- \frac{n}{p}}\leq Cs|\tilde\gamma - \frac{g(\xi,s)\tilde\gamma + \xi}{s}|^{1- \frac{n}{p}}\]
\[|T_4|\leq Cs (\|\tilde W_1\|_{C^{0,1 - \frac{n}{p}}(\tilde\Omega)} + \|\tilde W_2\|_{C^{0,1 - \frac{n}{p}}(\tilde\Omega)})| 1- \frac{g(\xi,s)}{s}|\leq Cs| 1- \frac{g(\xi,s)}{s}| \]
Substitute the identity for $T_1$ into equation(\ref{plus then minus}) and set $G(\xi,s) := T_2+ T_3+ T_4$ we get the desired result.$\square$



\noindent{\bf Proof of Theorem \ref{curl estimate}}\\
By \cite{sun}, any solution $u_l$ of $H_{W_l,q_l}u_l = 0$ satisfies the following identity

\begin{eqnarray}
\label{identity}
&&i\int_\Omega(W_1 - W_2)\cdot(u_1\nabla\bar u_2-\bar u_2\nabla u_1) + (W_1^2 - W_2^2)u_1\bar u_2 + (q_1 - q_2)u_1\bar u_2\\\nonumber &&= \int_{\partial\Omega} \bar u_2(\Lambda_{W_2,q_2}-\Lambda_{W_1,q_1})u_1
\end{eqnarray}
Combining the result of lemmas \ref{curl identity} and \ref{leading term estimate} into this equation we get that
\begin{eqnarray}
\label{fourier transform identity}
2s\int_{\Omega}\bar\mu\cdot(W_1 - W_2)e^{2i\xi\cdot x} &=& \int_{\partial\Omega}\bar u_2(\Lambda_{W_1,q_1} - \Lambda_{W_2,q_2})u_1\\\nonumber
&+& \int_\Omega(W_1^2 - W_2^2 + q_2 - q_1)\bar u_2u_1 + f_1(x,s)+ f_2(x,s) dx + G(\xi,s)  
\end{eqnarray}  
where $f_1$, $f_2$, and $G(\xi,s)$ are as in lemma \ref{curl identity} and \ref{leading term estimate}. Observe that if we set 
\[D := sup\{|x|\mid x\in\Omega\}\]
then $\|u_l\|_{H^1(\Omega)}\leq Ce^{Ds}$ for $s$ large. Furthermore, the $L^1(\Omega)$ norm of $|\bar u_2 u_1|$ is uniformly bounded independent of $s$ and $\xi$. Apply lemma \ref{leading term estimate} and lemma \ref{curl identity} we conclude the following estimate for all $\xi\in\R^n$, $s >\frac{1}{h_0}$ such that $|\xi_k| + |\xi_j|>0$ and $|\xi|\leq s$:
\begin{eqnarray*}
|\int_\Omega \bar\mu\cdot(W_1 - W_2) e^{2i\xi\cdot x}| 
&\leq& C(e^{Ds}\|\Lambda_{W_1,q_1} - \Lambda_{W_2,q_2}\|_{\frac{1}{2},\frac{-1}{2}} + |1 - \frac{g(\xi,s)}{s}| + s^{-\epsilon} \\ &&\ \ \ \ \ \ + |\tilde\gamma - \frac{g(\xi,s)\tilde\gamma + \xi}{s}|^{1-\frac{n}{p}})\\
&\leq& C(e^{Ds}\|\Lambda_{W_1,q_1} - \Lambda_{W_2,q_2}\|_{\frac{1}{2},\frac{-1}{2}} + (\frac{|\xi|}{s})^{1-n/p} + s^{-\epsilon}
\end{eqnarray*}

Here we used the inequality $\sqrt{1-t}\geq 1-t$ for $t\in[0,1]$. By the exact same method as above we can obtain the estimate for $-\mu$ in place of $\bar\mu$ to get 
\begin{eqnarray*}
|\int_\Omega -\mu\cdot(W_1 - W_2) e^{2i\xi\cdot x}| \leq C(e^{Ds}\|\Lambda_{W_1,q_1} - \Lambda_{W_2,q_2}\|_{\frac{1}{2},\frac{-1}{2}} + (\frac{|\xi|}{s})^{1-n/p} + s^{-\epsilon}
\end{eqnarray*}
Recall that $\mu = \gamma + i\tilde \gamma$ with $\tilde \gamma = \frac{\xi_j e_k - \xi_k e_j}{|\xi_je_k -\xi_ke_j|}$. Multiplying both estimates by $|\xi_je_k -\xi_ke_j|$ then add them together we get that:
\begin{eqnarray}
\label{not quite fourier}
|\int_\Omega e^{2i\xi\cdot x}\{\xi_j (W_1 - W_2)_k - \xi_k(W_1 - W_2)_j\}|\leq C|\xi|(e^{Ds}\|\Lambda_{W_1,q_1} - \Lambda_{W_2,q_2}\|_{\frac{1}{2},\frac{-1}{2}} + (\frac{|\xi|}{s})^{1-n/p} + s^{-\epsilon})
\end{eqnarray}
Recall that since $(W_1 - W_2) = 0$ on $\partial\Omega$, we can extend $(W_1 - W_2)$ to a $H^1(\R^n)$ vector field by defining it to be zero outside of $\Omega$ and we will refer to the extention as $I_\Omega(W_1 - W_2)$. With this extention we have that $d(I_\Omega(W_1 - W_2))$ is an $L^2(\R^n)$ function supported only in $\Omega$ and 
\[d(I_\Omega(W_1 -W_2)) = I_\Omega d(W_1 - W_2)\] 
as $L^2(\R^n)$ functions. Therefore, (\ref{not quite fourier}) implies that for all $\xi\in\R^n$ satisfying $|\xi_j| + |\xi_k| > 0$ and $s\geq \frac{1}{h_0}$ such that $|\xi|\leq s$, the following inequality holds for the Fourier transform of each component of $I_\Omega d(W_1-W_2)$
\begin{eqnarray*}
|{\cal F}(I_\Omega(\partial_j(W_1 - W_2)_k - \partial_k(W_1 - W_2)_j))(2\xi)|\leq C|\xi|(e^{Ds}\|\Lambda_{W_1,q_1} - \Lambda_{W_2,q_2}\|_{\frac{1}{2},\frac{-1}{2}} + (\frac{|\xi|}{s})^{1-n/p} + s^{-\epsilon})
\end{eqnarray*}

By the fact that both the right and left side are continuous, the estimate holds for all $\xi$ such that $|\xi|\leq s$. Since this is true for all components of the 2-form $I_\Omega d(W_1 - W_2)$ we have that for any $0<R\leq s$, $s\geq \frac{1}{h_0}$
\begin{eqnarray*}
\|I_\Omega d(W_1 - W_2)\|_{H^{-1}}^2  &=& \int_{|\xi|\leq R}\frac{|{\cal F}(I_\Omega d(W_1 - W_2))|^2}{{1+|\xi|^2}} d\xi+ \int_{|\xi|\geq R}\frac{|{\cal F}(I_\Omega d(W_1 - W_2))|^2}{{1+|\xi|^2}}d\xi\\
&\leq& CR^n(e^{2Ds}\|\Lambda_{W_1,q_1} - \Lambda_{W_2,q_2}\|_{\frac{1}{2},\frac{-1}{2}}^2 + (\frac{R^2}{s^2})^{1-n/p} + s^{-2\epsilon}) +\frac{C}{R^2}
\end{eqnarray*}
Here $\|\cdot\|_{H^{-1}}$ denotes the norm on $H^{-1}(\R^n)$, the dual space of $H^1(\R^n)$ and the last inequality comes from our a-priori assumption on the $W^{1,p}(\Omega)$ norm of $W_l$.
Recall however that the above statement is only valid for $s$ large enough to guarantee CGO solutions. Namely, the inequality holds only when $s\geq \frac{1}{h_0}$ where $h_0>0$ depends only on $\Omega$ and $M$. But since the estimate is trivial in the case when 
\[\|\Lambda_{W_1,q_1} - \Lambda_{W_2,q_2}\|_{\frac{1}{2},\frac{-1}{2}}\geq min\{1,e^{-2D/h_0}\}\] 
(just take the constant large enough and use the a-priori assumptions on $W_l$), we may assume without loss of generality that $\|\Lambda_{W_1,q_1} - \Lambda_{W_2,q_2}\|_{\frac{1}{2},\frac{-1}{2}}< min\{1,e^{-2D/h_0}\}$. With this assumption we may choose $R^n = s^{min\{\epsilon, 1-n/p\}}$ and $s = \frac{-log\|\Lambda_{W_1,q_1} - \Lambda_{W_2,q_2}\|}{2D}$ and we obtain that 
\[\|I_\Omega d(W_1 - W_2)\|_{H^{-1}}^2  \leq \frac{C}{|log\|\Lambda_{W_1,q_1} - \Lambda_{W_2,q_2}\|_{\frac{1}{2},\frac{-1}{2}}|^{\epsilon}} + C\|\Lambda_{W_1,q_1} - \Lambda_{W_2,q_2}\|_{\frac{1}{2},\frac{-1}{2}} \] 
for some $\epsilon > 0$, $C>0$ depending on $M$ and $\Omega$.$\square$\\\\
Note here that we have obtained an estimate for the $H^{-1}(\R^n)$ norm of $I_\Omega d(W_1 - W_2)$ and not just the $H^{-1}(\Omega)$ norm of $d(W_1 - W_2)$. For clarity, we will refer to the dual space of $H^1_0(\Omega)$ as $H^{-1}(\Omega)$ and the dual space to $H^1(\Omega)$ as $H^1(\Omega)^*$. In general, the $H^1(\Omega)^*$ norm is larger than the $H^{-1}(\Omega)$ norm. It is easily seen from the estimate which we derived in theorem \ref{curl estimate} and corollary \ref{reduction} that in both cases we have
\begin{eqnarray}
\label{H^1-star norm estimate}
\|d(W_1 - W_2)\|^2_{H^1(\Omega)^*} \leq \frac{C}{|log\|\Lambda_{W_1,q_1} - \Lambda_{W_2,q_2}\|_{\frac{1}{2},\frac{-1}{2}}|^{\epsilon}} + C\|\Lambda_{W_1,q_1} - \Lambda_{W_2,q_2}\|_{\frac{1}{2},\frac{-1}{2}}
\end{eqnarray}
This will be a key ingredient in proving the result of the next section.

\end{subsection}
\end{section}

\begin{section}{Estimates for the Electric Potential}
The goal of this section is to prove the following proposition:
\begin{theorem}
\label{potential estimate}
Let $\Omega$ be a bounded open subset of $\R^n$ with smooth boundary. For all $M>0$, there exists a $C> 0$, $\epsilon>0$ such that the estimate
\[\|q_1 - q_2\|_{H^{-1}(\R^n)}\leq C\lbrace \|\Lambda_{W_1,q_1}-\Lambda_{W_2,q_2}\|_{\frac{1}{2},\frac{-1}{2}}^{1/2} + |log \|\Lambda_{W_1,q_1}-\Lambda_{W_2,q_2}\|_{\frac{1}{2},\frac{-1}{2}}|^{-\epsilon}\rbrace\]

holds for all $W_1, W_2\in W^{2,\infty}(\Omega)$ and $q_1,q_2\in L^\infty(\Omega)$ satisfying $\|W_l\|_{ W^{2,\infty}(\Omega)}\leq M$, $\|q_l\|_{L^\infty(\Omega)}\leq M$ ($l = 1,2$) and $W_1\mid_D = W_2\mid_D$ on $\partial\Omega$. 

\end{theorem}
The proof of theorem \ref{potential estimate} involves using the stability result we already obtained for the magnetic field. In order to do this we first need to prove a lemma about the bounded invertibility of the operator $d$ on the set of 1-forms. This of course, is not true in general since the operator always has a non-trivial null space. But we will see that it is indeed injective provided that we quotient out the exact forms. We will employ the following notations. We consider $(\Omega,\partial\Omega)$ to be a Riemanian manifold with boundary and denote by $F^k(\Omega)$ to be the the set of k-forms on $\Omega$ and $W^{p,m}F^k(\Omega)$ to be its $W^{m,p}$ closure. Set
\[{\cal E}^k(\Omega) :=\{d\alpha\mid \alpha\in H^{1}_D F^{k-1}(\Omega)\},\ \ \ {\cal C}^k(\Omega) := \{d\alpha\mid \alpha\in H^{1}_N F^{k-1}(\Omega)\}\]
where $H^{1}_D F^{k}(\Omega)$ and $H^{1}_N F^{k}(\Omega)$ are the set of $H^1$ k-forms with homogenous Dirichlet and Neumann boundary trace, respectively. Furthermore we denote by ${\cal H}^k(\Omega)$ to be the $L^2$ closure of the space of harmonic k-forms. The corresponding subspaces in $W^{p,m}F^k(\Omega)$ are denoted by
\[W^{m,p}{\cal E}^k(\Omega) := {\cal E}^k(\Omega) \cap W^{m,p}F^k(\Omega),\ \ W^{m,p}{\cal C}^k(\Omega) := {\cal C}^k(\Omega) \cap W^{m,p}F^k(\Omega)\]
\[and\ \ \ W^{m,p}{\cal H}^k(\Omega) := {\cal H}^k(\Omega) \cap W^{m,p}F^k(\Omega)\]
We will identify the space of 1-forms $W^{m,p}F^1(\Omega)$ with the space of vector fields $W^{m,p}(\Omega;\R^n)$ and the space of 0-forms $W^{m,p}F^0(\Omega)$ with the space of functions $W^{m,p}(\Omega;\R)$.
\begin{lemma}
\label{bounded inverse of d}
Suppose $\Omega\subset \R^n$ is a bounded open set that is simply connected with connected boundary. For $p\geq 2$ and $m\in\N\backslash\{0\}$, define the set of Dirichlet 1- forms (or vector fields)
\[W^{m,p}_D(\Omega;\R^n) = \{W\in W^{m,p}(\Omega;\R^n)\mid W(x)\cdot\hat t = 0\ \forall \hat t\in T_x\partial\Omega,\ x\in\partial\Omega\}\] 
and 
\[X_0 = (W^{m,p}{\cal C}^1(\Omega) \oplus W^{m,p}{\cal H}^1(\Omega)) \cap W^{m,p}_D(\Omega;\R^n)\]
Then the differential $d: X_0 \to W^{m-1,p}{\cal E}^2(\Omega)$ has a bounded inverse. 
\end{lemma}
\noindent{\bf Proof}\\
Since $d: W^{m,p}F^k(\Omega) \to W^{m-1,p}F^{k+1}(\Omega)$ is a continuous linear operator we will use uniform boundedness principle to assert the existence of a bounded inverse. Since by \cite{schwarz}, $W^{m,p}{\cal C}^1(\Omega)$, $W^{m,p}{\cal H}^1(\Omega)$, and $W^{m,p}_D(\Omega;\R^n)$ are all closed subspaces of $W^{m,p}F^1(\Omega)$, it is clear that $X_0$ is a closed subspace of $W^{m,p}F^1(\Omega)$. Furthermore, $W^{m-1,p}{\cal E}^2(\Omega)$ is a closed subspace of $W^{m-1,p}F^2(\Omega)$ so it suffices to check the bijectivity of $d: X_0\to W^{m-1,p}{\cal E}^2(\Omega)$ for the uniform boundedness principle to apply.\\
To see injectivity, suppose $dW = 0$ for some $W\in X_0$. By the fact that cohomology of $\Omega$ is assumed to be trivial, this means that $W = d\alpha$ for some $\alpha\in W^{m+1,p}(\Omega;\R)$. Since $W$ has no tangential component and $\partial\Omega$ is connected this means that $\alpha$ is a constant function along the boundary and we may take $\alpha\in W^{m+1,p}_D(\Omega)$. By definition this means that 
\[W\in W^{m,p}{\cal E}^1(\Omega) := {\cal E}^1(\Omega)\cap W^{m,p}F^1(\Omega)\]
However, $W\in X_0\subset W^{m,p}{\cal C}^1(\Omega) \oplus W^{m,p}{\cal H}^1(\Omega)$ which is $L^2$ perpendicular to $W^{m,p}{\cal E}^1(\Omega)$ and thus $W=0$. So this establishes injectivity.\\
To show surjectivity, let $\omega \in  W^{m-1,p}{\cal E}^2(\Omega)$. Then by definition, $\omega = dW$ for some $W\in W^{1,2}_D(\Omega;\R^n)$. By the Hodge decomposition, we can decompose $\omega = dW_\omega + \delta\beta_\omega + \kappa_\omega$ with $W_\omega \in W^{m,p}_D(\Omega;\R^n)$. Subtract the two expressions we have for $\omega$ and use $L^2$ orthogonality again we see that 
\[dW_\omega = dW = \omega\]
We now have that $\omega$ is in the range of the $d$ operator acting on $W_D^{m,p}(\Omega;\R^n)$ but it is not clear that $W_\omega\in X_0$. This can be remedied by applying the Hodge decomposition to $W_\omega$ to produce 
\begin{eqnarray}
\label{decomposition of W}
W_\omega - d\alpha = \delta\eta + \kappa\ \ \ with\ \ \ \alpha\in W^{m+1,p}_D(\Omega;\R),\ \ \delta\eta\in W^{m,p}{\cal C}^1(\Omega),\ \ \kappa\in W^{m,p}{\cal H}^1(\Omega)
\end{eqnarray}
Since $\alpha\in W^{m+1,p}_D(\Omega;\R)$, $d\alpha$ has no tangential component along the boundary. This with \\${W_\omega\in W_D^{m,p}(\Omega;\R^n)}$ shows that the LHS of (\ref{decomposition of W}) is in $W_D^{m,p}(\Omega;\R^n)$. The RHS of (\ref{decomposition of W}) is clearly in \\$ W^{m,p}{\cal C}^1(\Omega)\oplus W^{m,p}{\cal H}^1(\Omega)$ so we conclude that $W_\omega - d\alpha\in X_0$. We now have 
\[\omega = dW_\omega = d(W_\Omega - d\alpha)\ \  with\ \ W_\Omega - d\alpha\in X_0\] 
So surjectivity is established and uniform boundedness principle applies to give a bounded inverse.$\square$\\

We would like to apply this lemma directly to the vector field $(W_1 - W_2)$ which may not be in $X_0$. So first we do the following manipulation. Pick $p> {p_0} > n$ and apply the Hodge decomposition to $W_1 - W_2$ in the space $W^{2,{p_0}}(\Omega;\R^n)$ to get $W_1 - W_ 2 = \delta\beta + d\alpha + \kappa$. By lemma 2.4.11 of \cite{schwarz}, $\alpha \in W^{3,{p_0}}(\Omega)\cap H^1_0(\Omega)$ with norm $\|\alpha\|_{W^{3,{p_0}}(\Omega)}\leq C\|W_1 - W_2\|_{W^{2,{p_0}}(\Omega;\R^n)}\leq 2CM$ where the constant $C$ depends on $\Omega$ and ${p_0}$ only.  Define $W_1' := W_1 -\frac{d\alpha}{2}$ and $W_2' := W_2 + \frac{d\alpha}{2}$ so that $W_1' - W_2' \in W^{2,{p_0}}{\cal C}^1(\Omega) \oplus W^{2,{p_0}}{\cal H}^1(\Omega)$. Since we have already assumed that $W_1 - W_2$ has no tangential component at the boundary and $\alpha\in H^1_0(\Omega)$, we conclude that $W_1' - W_2'$ has no tangential component. This means that 
\[W_1' - W_2'\in (W^{1,{p_0}}{\cal C}^1(\Omega) \oplus W^{1,{p_0}}{\cal H}^1(\Omega))\cap H^1_D(\Omega;\R^n)\] 
so by lemma \ref{bounded inverse of d} 
\begin{eqnarray}
\label{curl trick}
\|W_1' - W_2'\|_{W^{1,{p_0}}(\Omega)}\leq C\|d(W_1' - W_2')\|_{L^{p_0}(\Omega)} = C\|d(W_1 - W_2)\|_{L^{p_0}(\Omega)}
\end{eqnarray}
Recall that due to gauge invariance, we have $\Lambda_{W_l,q_l} = \Lambda_{W_l',q_l}$.\\
Choose a bounded open set $\tilde\Omega$ such that $\Omega\subset\subset\tilde\Omega$. Let $E$ be the extension operator mapping $W^{1,{p_0}}(\Omega) \to W^{1,{p_0}}_c(\tilde\Omega)$ and $W^{2,{p_0}}(\Omega)\to W^{2,{p_0}}_c(\tilde\Omega)$ such that it is bounded in both norms (see \cite{evans}). Denote by $\tilde W_l' := E W_l'$ for $(l = 1,2)$. Note that the extention described here is different from the one we used in the previous section. For one thing we no longer have $\tilde W_1' = \tilde W_2'$ in $\tilde\Omega\backslash\Omega$. However, we still have that for each $M>0$ there exists an $M'>0$ and $R>0$ such that if $\|W\|_{W^{2,{p_0}}(\Omega)}\leq M$ then $\tilde W'\in W^{2,p_0}_c(\tilde\Omega)$ with $\|\tilde W'\|_{W^{2,p_0}(\tilde\Omega)}\leq M'$.\\
In this setting proposition \ref{existence of CGO} applied to $\tilde W_l'$ gives a $C>0$, $(1-\frac{n}{p_0})\sigma_0>\epsilon >0$, and $h_0>0$ such that if $\zeta_l\in\C^n$ satisfies $\zeta_l\cdot\zeta_l = 0$, $|\zeta_l|\geq 1/h_0$ then for all $W_l,q_l$ $(l = 1,2)$ satisfying $\|W_l\|_{W^{2,p_0}(\Omega;\R^n)}\leq M$ and $\|q_l\|_{L^\infty}\leq M$ there exists solutions to $H_{\tilde W'_l,q_l}u_l = 0$ of the form
\[u_l(x,\zeta_l) = e^{i\zeta_l\cdot x}(e^{i\chi_{|\zeta_l|}\phi^{\sharp}_l} + r_l(x,\zeta_l))\ \ \ \|r_l(\cdot,\zeta)\|_{H^t_\delta}\leq C|\zeta_l|^{t-\epsilon}\ \ 
\ \ t \in [0,2]\ \ -1<\delta<0\]


We emphasize again that the constants $C>0$, $\epsilon >0$ and $h_0>0$ depends as usual on the parameters and the a-priori bound $M$ but not on the choice of $(W_l,q_l)$ that satisfies the a-priori assumption.\\
\noindent{\bf Proof of Theorem \ref{potential estimate}}\\
We start again with identity (\ref{identity}) except this time we will isolate the electric potential term on the LHS. Similar calculation to those done in lemma \ref{leading term estimate} shows that 
\[\|u_1 \bar u_2\|_{L^1(\Omega)}\leq C,\ \ \ \|u_1\nabla\bar u_2\|_{L^1(\Omega)} + \|\bar u_2\nabla u_1\|_{L^1(\Omega)} \leq Cs\]
Use these estimates and the fact that $\|r_l(x,s)\|_{L^2(\Omega)}\leq Cs^{-\epsilon}$ identity (\ref{identity}), becomes
\begin{eqnarray*}
|\int_\Omega (q_1 - q_2)e^{i\xi\cdot x}e^{i\phi_1^\sharp - i\overline{\phi_2^\sharp}}|\leq
 C(e^{2Ds}\|\Lambda_{W_1,q_1}-\Lambda_{W_2,q_2}\|_{\frac{1}{2},\frac{-1}{2}} + \|W_1'-W_2'\|_{L^\infty(\Omega)}s +s^{-\epsilon})
\end{eqnarray*}
Here again we use the fact that when $s$ is large the CGO solutions satisfy $\|u_l(x,\zeta)\|_{H^1(\Omega)}\leq Ce^{sD}$ where 
$D := sup\{|x|\mid x\in\Omega\}$.
Apply Morrey's inequality to $\|W_1' - W_2'\|_{L^\infty(\Omega)}$ then use (\ref{curl trick}) we get that 
\begin{eqnarray*}
|\int_\Omega (q_1 - q_2)e^{i\xi\cdot x}e^{i\phi_1^\sharp - i\overline{\phi_2^\sharp}}|&\leq& C(e^{2Ds}\|\Lambda_{W_1,q_1}-\Lambda_{W_2,q_2}\|_{\frac{1}{2},\frac{-1}{2}} + \|d(W_1-W_2)\|_{L^p(\Omega)}s + s^{-\epsilon})
\end{eqnarray*}
We want a Fourier transform to appear on the LHS. Therefore we replace $\int_\Omega (q_1 - q_2)e^{i\xi\cdot x}e^{i\phi_1^\sharp - i\overline{\phi_2^\sharp}}$ by $\int_\Omega (q_1 - q_2)e^{i\xi\cdot x}$ to obtain
\begin{eqnarray}
\label{fourier transform of q}
\nonumber|{\cal F}(q_1 - q_2)(2\xi)|&\leq& C(e^{2Ds}\|\Lambda_{W_1,q_1}-\Lambda_{W_2,q_2}\|_{\frac{1}{2},\frac{-1}{2}} + \|d(W_1-W_2)\|_{L^{p_0}(\Omega)}s + s^{-\epsilon}\\&&\ \ \ \ \ \ \ \ \ \ \ \ \ \ \ \ \ \ \ \ \ \ \ \ \ \ \ \ \ \ \ \ \ \ \ \ \ \ \ \ \ \ \ \ \ \ \ \ \ \ \ \ \ \ \ \ \ \ \ \ \ \ \ \ \  + \|\phi_1^\sharp - \overline{\phi_2^\sharp}\|_{L^\infty(\Omega)})
\end{eqnarray}
with all constants depending only on $M$.\\
Now we would like to estimate the last term by $\|W_1' - W_2'\|_{W^{1,p}}$. This can be done by writing
\begin{eqnarray}
\label{triangle}
\|\phi_1^\sharp - \overline{\phi_2^\sharp}\|_{L^\infty(\Omega)}\leq \|\phi_1^\sharp - {\psi_1^\sharp}\|_{L^\infty(\Omega)} + \|\overline{\phi_2^\sharp} - \overline{\psi_2^\sharp}\|_{L^\infty(\Omega)} + \|\psi_1^\sharp - \overline{\psi_2^\sharp}\|_{L^\infty(\Omega)}
\end{eqnarray} 
By our extention, $\tilde W_l'\in W^{2,p_0}_c(\tilde\Omega)\hookrightarrow C_c^{1}(\tilde\Omega)$. Therefore lemma \ref{uniform estimate} shows that the first two terms are bounded by $C\sqrt{(1-\sqrt{1-|\xi|^2/s^2})}$ with $C$ depending only on the a-priori bound $M$. Due to lemma \ref{L-infinity bound} the last term is bounded by 
\[\|\psi_1^\sharp - \overline{\psi_2^\sharp}\|_{L^\infty(\Omega)}\leq C\|\tilde W'^\sharp_1 - \tilde W'^\sharp_2\|_{L^{\infty}(\tilde\Omega)} \leq C\|\tilde W'_1 - \tilde W'_2\|_{L^{\infty}(\tilde\Omega)}+ \| \tilde W_2'^\flat\|_{L^\infty(\tilde\Omega)}+ \| \tilde W_1'^\flat\|_{L^\infty(\tilde\Omega)}\] 

Substitute this into (\ref{triangle}) and use Morrey inequality we get that
\begin{eqnarray*}
\|\phi_1^\sharp - \overline{\phi_2^\sharp}\|_{L^\infty(\Omega)}&\leq& C\sqrt{1-\sqrt{1-|\xi|^2/s^2}} + C\|\tilde W'_1 - \tilde W'_2\|_{W^{1,p_0}(\tilde\Omega)}+ \| \tilde W_2'^\flat\|_{L^\infty(\tilde\Omega)}+ \|\tilde W_1'^\flat\|_{L^\infty(\tilde\Omega)}\\
&\leq& C\sqrt{1-\sqrt{1-|\xi|^2/s^2}} + C\| W'_1 -  W'_2\|_{W^{1,p_0}(\Omega)}+ + \| \tilde W_2'^\flat\|_{L^\infty(\tilde\Omega)}+ \| \tilde W_1'^\flat\|_{L^\infty(\tilde\Omega)}\\
&\leq&C\sqrt{1-\sqrt{1-|\xi|^2/s^2}} + C\| d(W_1 -  W_2)\|_{L^{p_0}(\Omega)}+ + \| \tilde W_2'^\flat\|_{L^\infty(\tilde\Omega)}+ \| \tilde W_1'^\flat\|_{L^\infty(\tilde\Omega)}
\end{eqnarray*}
Use the fact that $\tilde W_l'\in W^{2,p_0}_c(\tilde\Omega)\hookrightarrow C_c^{1}(\tilde\Omega)$ again we see that $ \| \tilde W_l'^\flat\|_{L^\infty(\tilde\Omega)}$ vanishes like $Cs^{-\sigma_0}$ with $C$ depending only on the a-priori bound $M$. Substitute this into the inequality (\ref{fourier transform of q}) we get
\begin{eqnarray*}
|{\cal F}(q_1 - q_2)(2\xi)|\leq C(e^{2Ds}\|\Lambda_{W_1,q_1}-\Lambda_{W_2,q_2}\|_{\frac{1}{2},\frac{-1}{2}} + \|d(W_1-W_2)\|_{L^{p_0}(\Omega)}s + s^{-\epsilon})
\end{eqnarray*}
Recall that $p_0$ and $p$ satisfy the condition $n <p_0<p$ so we can interpolate 
\[\|d(W_1-W_2)\|_{L^{p_0}(\Omega)}\leq \|d(W_1-W_2)\|^{1-\lambda}_{L^{p}(\Omega)}\|d(W_1-W_2)\|^{\lambda}_{L^{2}(\Omega)},\ \lambda = \frac{p_0^{-1}-p^{-1}}{2^{-1}-p^{-1}}\] 

and use our a-priori assumption about the $W^{2,p}(\Omega)$ norm of $W_l$ we get that
\begin{eqnarray}
\label{fourier transform of q II}
|{\cal F}(q_1 - q_2)(2\xi)|\leq C(e^{2Ds}\|\Lambda_{W_1,q_1}-\Lambda_{W_2,q_2}\|_{\frac{1}{2},\frac{-1}{2}} + \|d(W_1-W_2)\|^\lambda_{L^{2}(\Omega)}s + s^{-\epsilon})
\end{eqnarray}
For any $f\in H^1(\Omega)$, $\|f\|_{L^2(\Omega)}^2\leq \|f\|_{H^1(\Omega)}\|f\|_{H^{1}(\Omega)^*}$. Apply this to each component of the 2-form $d(W_1 - W_2)$ and use our a-priori bound on $\|W_l\|_{W^{2,p}(\Omega)}$, (\ref{fourier transform of q II}) becomes
\[|{\cal F}(q_1 - q_2)(2\xi)|\leq C(e^{2Ds}\|\Lambda_{W_1,q_1}-\Lambda_{W_2,q_2}\|_{\frac{1}{2},\frac{-1}{2}} + \|d(W_1-W_2)\|_{H^1(\Omega)^*}^{\lambda/2}s + s^{-\epsilon})\] 
So by estimate (\ref{H^1-star norm estimate}) we have
\begin{eqnarray}
\label{last estimate on fourier transform}
\nonumber|{\cal F}(q_1 - q_2)(2\xi)|&\leq&C(e^{2Ds}\|\Lambda_{W_1,q_1}-\Lambda_{W_2,q_2}\|_{\frac{1}{2},\frac{-1}{2}} +\\ &&\ \ \{|log\|\Lambda_{W_1,q_1} - \Lambda_{W_2,q_2}\|_{\frac{1}{2},\frac{-1}{2}}|^{-\epsilon} + \|\Lambda_{W_1,q_1} - \Lambda_{W_2,q_2}\|_{\frac{1}{2},\frac{-1}{2}}\}^{\lambda/4} s + s^{-\epsilon})
\end{eqnarray}
and this is valid for all $s \geq \frac{1}{h_0}$ and $s \geq |\xi|$. Computing directly the $H^{-1}(\R^n)$ norm of $q_1 - q_2$ by using Fourier transform, (\ref{last estimate on fourier transform}) implies that for $R>0$, $s \geq max\{R,\frac{1}{h_0}\}$
\begin{eqnarray*}
\|q_1 - q_2\|^2_{H^{-1}}\leq R^n C&(&e^{2Ds}\|\Lambda_{W_1,q_1}-\Lambda_{W_2,q_2}\|_{\frac{1}{2},\frac{-1}{2}} +\\&&\ \   \{|log\|\Lambda_{W_1,q_1} - \Lambda_{W_2,q_2}\|_{\frac{1}{2},\frac{-1}{2}}|^{-\epsilon} + \|\Lambda_{W_1,q_1} - \Lambda_{W_2,q_2}\|_{\frac{1}{2},\frac{-1}{2}}\}^{\lambda/4} s + s^{-\epsilon})+ \frac{C}{R^2} 
\end{eqnarray*}
Clearly, we can take $\epsilon < 1$ without loss of generality. Furthermore, as in the previous section we may assume that $\|\Lambda_{W_1,q_1}-\Lambda_{W_2,q_2}\|_{\frac{1}{2},\frac{-1}{2}}\leq min\{e^{-(\frac{4D}{h_0})^{8/(\lambda\epsilon)}},1\}$ since the quantity $e^{-(\frac{4D}{h_0})^{8/(\lambda\epsilon)}}$ depends only on the a-priori bound $M$ but not the chosen electric and magnetic potentials. Under these assumptions we may take $R^n = s^{\epsilon/2}$ and $s = \frac{1}{4D}|log\|\Lambda_{W_1,q_1}-\Lambda_{W_2,q_2}\|_{\frac{1}{2},\frac{-1}{2}}|^{\epsilon/8}\geq \frac{1}{h_0}$ and simple arithmetics establishes the lemma.$\square$

\end{section}
\end{section}
\begin{section}*{PART III- Partial Data Estimate}
\begin{section}{Introduction}
We now move on to the third part of the paper which addresses the stability problem when DN map is known only on part of the boundary. The estimates will be of $log\ log$-type and will be weaker than the one established in the full data case. For each $\tilde\gamma\in S^{n-1}$ and $\epsilon_0 >0$ we define the front and back of the boundary $\partial\Omega$ with respect to $\tilde\gamma$ by
\[\partial\Omega_{+,\epsilon_0}(\tilde\gamma)= \{x\in\partial\Omega \mid \tilde\gamma\cdot\nu(x) >\epsilon_0\}\ \ \ \partial\Omega_{-,\epsilon_0}(\tilde\gamma)= \partial\Omega\backslash\overline{\partial\Omega_{+,\epsilon_0}(\tilde\gamma)}\]
Using this notation with fixed $\epsilon_0>0$, we define for any magnetic potential $W$ and electric potential $q$ the partial DN map $\Lambda_{W,q}' :H^{3/2}(\partial\Omega) \to H^{1/2}(\partial\Omega_{-,2\epsilon_0}(e_n))$ by
\[\Lambda_{W,q}'f = (\Lambda_{W,q}f)\mid_{\partial\Omega_{-,2\epsilon_0}(e_n)}\ \ \ \ \forall f\in H^{3/2}(\partial\Omega)\]
with the associated operator norm
\[\|\Lambda_{W,q}'\|' := sup\{\|\Lambda_{W,q}'f\|_{H^{1/2}(\partial\Omega_{-,2\epsilon_0}(e_n))} \mid \|f\|_{H^{3/2}(\partial\Omega)}=1\}\]
Given two sets of electric and magnetic potentials $(W_1,q_1)$ and $(W_2,q_2)$ we wish to estimate their difference in terms of $\|\Lambda_{W_1,q_1}' - \Lambda_{W_2,q_2}'\|'$. More precisely,
\begin{theorem}
\label{partial data stability for magnetic field}
Let $\Omega$ be a bounded open subset of $\R^n$ with smooth boundary. For all $M>0$, there exists a $C> 0$, $\lambda>0$ such that the estimate
\[\|I_{\Omega}d(W_1 - W_2)\|_{H^{-1}(\R^n)}^2\leq C\lbrace \|\Lambda_{W_1,q_1}'-\Lambda_{W_2,q_2}'\|'^\lambda + |log |log \|\Lambda_{W_1,q_1}'-\Lambda_{W_2,q_2}'\|'||^{-2}\rbrace\]
holds for all $W_1, W_2\in W^{1,\infty}(\Omega)$ and $q_1,q_2\in L^\infty(\Omega)$ satisfying $\|W_l\|_{ W^{1,\infty}(\Omega)}\leq M$, $\|q_l\|_{L^\infty(\Omega)}\leq M$ ($l = 1,2$) and $W_1 = W_2$ on $\partial\Omega$. Here $I_\Omega(x)$ is the indicator function of $\Omega$.
\end{theorem}

\begin{theorem}
\label{partial data stability for electric potential}
Let $\Omega$ be a bounded open subset of $\R^n$ with smooth boundary. For all $M>0$, there exists a $C> 0$, $\lambda>0$ such that the estimate
\[\|I_\Omega(q_1 -q_2)\|_{H^{-1}(\R^n)}^2\leq C\lbrace \|\Lambda_{W_1,q_1}'-\Lambda_{W_2,q_2}'\|'^\lambda + |log |log \|\Lambda_{W_1,q_1}'-\Lambda_{W_2,q_2}'\|'||^{-\lambda}\rbrace\]
holds for all $W_1, W_2\in W^{1,\infty}(\Omega)$ and $q_1,q_2\in L^\infty(\Omega)$ satisfying $\|W_l\|_{ W^{1,\infty}(\Omega)}\leq M$, $\|q_l\|_{L^\infty(\Omega)}\leq M$ ($l = 1,2$) and $W_1 = W_2$ on $\partial\Omega$. Here $I_\Omega(x)$ is the indicator function of $\Omega$.
\end{theorem}
We conclude the introduction with a geometric observation which will be useful later. If we denote by $N(\delta) =\{\tilde \gamma\in S^{n-1}\mid d_{S^{n-1}}(\tilde\gamma,e_n) \leq \delta\}$ to be a $\delta$-neighbourhood around $e_n$ on $S^{n-1}$ then for each $\epsilon_0 >0$ there exists a $\delta >0$ such that 
\[\partial\Omega_{-,\epsilon_0}(\tilde\gamma)\subset\partial\Omega_{-,2\epsilon_0}(e_n),\ \ \ \forall\tilde\gamma\in N(\delta)\]
A simple geometric argument shows that for each $\delta>0$ there exists a $r_1 >0$ such that if $a\in\R,b>0$ are real numbers satisfying $|a|\leq r_1 b$  then $\frac{be_n -ae_j}{|be_n - ae_j|}\in N(\delta)\cap span\{e_n,e_j\}$. With this fact as motivation we define for each $j\in\{1,..,n-1\}$ the $j$-wedge by
\[E_j := \{\xi\in \R^n\mid |\xi_n|\leq r_1 \xi_j\}\]
Our proof of stability will be in two steps. First we will prove stability of the Fourier transform only in the $j$-wedge. We will then use a stable analytic continuation result to obtain stability of the Fourier transform in a ball of arbitrary radius. With this strategy in mind we will proceed with the next section on some necessary estimates.
\end{section}
\begin{section}{Preliminary Estimates}
In this section we state without proof two estimates which will be useful in deriving the partial data stability estimate. The proofs are given in the references. The first is a Caleman estimate which will allow us to bound the information we don't have on $\partial\Omega_{+,\epsilon_0}$ by the information we do have on $\partial\Omega_{-,\epsilon_0}$. We begin by first defining the notion of a limiting Carleman weight. 
\begin{definition}
\label{limiting weight}
A real smooth function $\phi$ on $\tilde\Omega$ is said to be a limiting Carleman weight if it has non-vanishing gradient on $\tilde\Omega$ and satisfies
\[< \phi''\nabla\phi,\nabla\phi> +  <\phi''\xi,\xi> =0\ \ \ \ when\ \ \xi^2 = (\nabla\phi)^2\ \ \ and\ \ \ \nabla\phi\cdot\xi = 0\]
\end{definition}
We now state the estimate we will use.
\begin{proposition}\cite{FKSU}
\label{carleman}
Let $\tilde\Omega$, $\Omega$ be bounded open subsets in $R^n$ such that $\Omega\subset\subset\tilde\Omega$ and $\phi$ be a $C^\infty$ limiting Carleman weight on $\tilde\Omega$. For every $M >0$, there exists a $C >0$, $h_0 >0$ such that for all $C^1$ vector fields $A$ on $\bar\Omega$ and $q\in L^\infty(\Omega)$ satisfying $max(\|A\|_{C^1(\Omega)}, \|q\|_{L^\infty(\Omega)})\leq M$, the following estimate holds  
\begin{eqnarray}
\label{carleman estimate}
-h(\partial_\nu \phi e^{\frac{\phi}{h}}\partial_\nu u\mid e^{\frac{\phi}{h}}\partial_\nu u)_{\partial\Omega_-} + \|e^{\frac{\phi}{h}}h\nabla u\|^2
\leq C h^2 \|e^{\frac{\phi}{h}} H_{A,q}u\|^2 + h(\partial_\nu\phi e^{\frac{\phi}{h}}\partial_\nu u\mid e^{\frac{\phi}{h}}\partial_\nu u)_{\partial\Omega_+}
\end{eqnarray}
for all $0<h<h_0$ and $u\in H^2(\Omega)\cap H^1_0(\Omega)$.
\end{proposition}
We will apply this proposition in the case when the limiting Carleman weight is $x\cdot\tilde\gamma$ with $\tilde\gamma\in N(\delta)$. The result will be the type of estimates on the difference of the Fourier transform that we have been seeing in part II. However, this time the estimate is only valid on a small wedge in phase space instead of an entire ball. To remedy this problem we follow the idea of Heck-Wang in \cite{hw} and use a stable dependence result for analytic continuation established by Vessella to extend the estimate on the wedge to an estimate on the ball.

\begin{proposition}\cite{vessella}
\label{analytic cont}
Let $B_2$ denote the open ball of radius $2$ centred around the origin and let $E$ be an open subset set of $B_1$. If $f$ is an analytic function with
\[\|\partial^\alpha f\|_{L^\infty(B_2)}\leq \frac{M\alpha !}{\rho^{|\alpha|}},\ \ \ \ \forall\alpha\in\N^n\]
for some $M,a>0$ then
\[\|f\|_{L^\infty(B_1)}\leq (2M)^{1- \theta(|E|/|D|)}(\|f\|_{L^\infty(E)})^{\theta(|E|/|D|)}\]
where $\theta \in (0,1)$ depends on $\rho$ and $n$.
\end{proposition}
We will see later that the neccesity of this estimate is what contributes to the $log\ log$ rate of convergence for the partial data as oppose to the $log$-type stability we have for the full data problem. I would be interesting to see whether one can refine such an estimate and consequently derive a $log$-type estimate for the partial data problem.

\end{section}

\begin{section}{Partial Data Estimate for the Magnetic Field}

We prove theorem \ref{partial data stability for magnetic field} in this section. We begin again with a discussion of (CGO) solutions. Note that this time we are only allowed to consider phase variables $\xi$ in the wedge $E_j$. This is due to the fact that we need $\tilde\gamma\in N(\delta)$ for the estimate in proposition \ref{carleman} to apply.\\
For all $\xi\in E_j$ $(j=1,..n-1)$ such that $\xi_j\not= 0$, define $\tilde \gamma := \frac{\xi_je_n - \xi_ne_j}{|\xi_je_n - \xi_ne_j|}\in N(\delta)$. Here $E_j$ and $\delta$ are chosen as in the introduction of Part II.\\
Choose $\gamma\perp\tilde\gamma$, $\gamma\perp\xi$, $ |\gamma| = 1$  and define $\mu = \gamma+i\tilde\gamma$. For each $s\geq|\xi|$ let
\[\zeta_1(s) = \zeta_1 = -is\tilde\gamma + g(\xi,s)\gamma + \xi\]
\[\zeta_2(s) = \zeta_2 = is\tilde\gamma + g(\xi,s)\gamma -\xi\]
where $g(\xi,s) = \sqrt{s^2 -|\xi|^2}$. Note that $\zeta_l(s)$ can be written as $\zeta_l(s) = (-1)^{l}is\tilde \gamma + s\gamma_l'$ $(l = 1,2)$ provided we take 
\[\gamma'_1 = \frac{g(\xi,s)\gamma + \xi}{s},\ \ \ \gamma'_2 = \frac{g(\xi,s)\gamma - \xi}{s}\]
Let $\sigma_0>0$ and $\theta>0$ be chosen so that $0<\sigma_0 <\sigma_0+\theta<\frac{1}{4n+6}$. Given any $W_1, W_2 \in W^{1,\infty}(\Omega)$ such that $W_1 = W_2$ on $\partial\Omega$ and $\|W_l\|_{W^{1,\infty}(\Omega)}\leq M$,  extend them to compactly supported vector fields $\tilde W_1, \tilde W_2\in W^{1,p}_c(\tilde\Omega;\R^n)$ as described in lemma \ref{divergence extention}. Therefore, there exists an $M'>0$ and $R>0$ such that $\tilde W_l \in {\cal W}(M',R)$ whenever $\|W_l\|_{W^{1,\infty}(\Omega)}\leq M$. Proposition \ref{existence of CGO} then gives a $C>0$, $h_0 >0$, and $ 0<\epsilon< \sigma_0(1-\frac{n}{p})$ such that for all $|\zeta_l|\geq\frac{1}{h_0}$ there exists solutions $u_l(x,\zeta_l)$ to $H_{\tilde W_l,q_l}u_l = 0$ in $\R^n$ that are of the form:
\[u_l(x,\zeta_l) = e^{i\zeta_l\cdot x}(e^{i\chi_{|\zeta_l|}\phi^{\sharp}_l} + r_l(x,\zeta_l))\ \ \ \|r_l(\cdot,\zeta_l)\|_{H^t_\delta}\leq C|\zeta_l|^{t-\epsilon}\ 
\ \ \ t\in[0,2]\ \ \ -1<\delta<0\]
\[\phi^\sharp_l(x) = -\frac{\gamma_l' +(-1)^li\tilde\gamma}{2\pi}\cdot\int_{\R^2}\frac{\tilde W_l^\sharp(x -\gamma_l'y_1 - (-1)^l\tilde\gamma y_2)}{y_1 + iy_2}\]
where $\chi_{|\zeta_l|} = \chi(x/|\zeta_l|^{\theta})$ and $\tilde W_l^\sharp = \int |\zeta_l|^{n\sigma_0}\chi(y|\zeta_l|^{\sigma_0}) \tilde W_l(x - y) dy$


Note here that $C$, $h_0$, and $\epsilon$ depend on dimension, $\Omega$, $\theta$, $\sigma_0$ and $M$ only but not on the choice of $(W_l,q_l)$ as long as $\|W_l\|_{W^{1,\infty}(\Omega)}\leq M$ and $\|q_l\|_{L^\infty}\leq M$. All constants in this section will have only the aforementioned dependence.

\begin{subsection}{Stability of Fourier Transform in a Cone}
In this section we derive an estimate for the Fourier transform ${\cal F}(\partial_j(\tilde W_1 - \tilde W_2)_n - \partial_n(\tilde W_1 - \tilde W_2)_j)(\xi)$ in the wedge $\xi\in E_j$. The main proposition is
\begin{lemma}
\label{wedge estimate}
The following estimate holds for all $\xi\in E_j$ such that $|\xi|\leq s$
\begin{eqnarray*}
&&|\int_\Omega e^{2i\xi\cdot x}\{\xi_j( W_1 - W_2)_n - \xi_n(W_1 - W_2)_j\}dx|\\
&&\leq  C|\xi|(e^{Ds}\|\Lambda'_{W_1,q_1} - \Lambda'_{W_2,q_2}\|' +(\sqrt{1-\sqrt{1-|\xi|^2/s^2}})^{1-\frac{n}{p}} + s^{-3\epsilon/4})
\end{eqnarray*}
\end{lemma}
\noindent{\bf Proof}\\
For each $\xi\in E_j$ let $\tilde\gamma := \frac{\xi_je_n - \xi_ne_j}{|\xi_je_n - \xi_ne_j|}\in N(\delta)$ and $u_1,u_2$ be the CGO solutions considered at the beginning of the section. We start with equation (\ref{fourier transform identity})
\begin{eqnarray}
\label{fourier transform identity II}
2s\int_{\Omega}\bar\mu\cdot(W_1 - W_2)e^{2i\xi\cdot x} &=& \int_{\partial\Omega}\bar u_2(\Lambda_{W_1,q_1} - \Lambda_{W_2,q_2})u_1\\\nonumber
&+& \int_\Omega(W_1^2 - W_2^2 + q_2 - q_1)\bar u_2u_1 + f_1(x,s)+ f_2(x,s) dx + G(\xi,s)  
\end{eqnarray}  
where $f_1$, $f_2$, and $G(\xi,s)$ are as in lemma \ref{curl identity} and \ref{leading term estimate}. From the estimates in these two lemmas we deduce that 
\begin{eqnarray}
\label{tail estimate}
\nonumber|\int_\Omega(W_1^2 - W_2^2 + q_2 - q_1)\bar u_2u_1 + f_1(x,s)&+& f_2(x,s) dx + G(\xi,s)|\\ \leq&& s(s^{-\epsilon} + (\sqrt{1-\sqrt{1-|\xi|^2/s^2}})^{1-\frac{n}{p}})
\end{eqnarray}
Denote by $v_l$ to be the solution of $H_{W_l,q_l}v_l = 0$ in $\Omega$ and $v_l = u_1$ on $\partial\Omega$. Then the first term of (\ref{fourier transform identity}) can be written as
\begin{eqnarray}
\label{the main term}
\nonumber |\int_{\partial\Omega}\bar u_2(\Lambda_{W_1,q_1} - \Lambda_{W_2,q_2})u_1| &=& |\int_{\partial\Omega_{-,\epsilon_0}(\tilde\gamma)}\bar u_2(\Lambda_{W_1,q_1}' - \Lambda_{W_2,q_2}')u_1+ \int_{\partial\Omega_{+,\epsilon_0}(\tilde\gamma)}\bar u_2\partial_\nu v|\\
&\leq& e^{2D}\|\Lambda_{W_1,q_1}' - \Lambda_{W_2,q_2}'\|' + |\int_{\partial\Omega_{+,\epsilon_0}(\tilde\gamma)}\bar u_2\partial_\nu v|
\end{eqnarray}
where $D:= sup\{|x|\mid x\in\Omega\}$ and $v := v_1 - v_2$. The particular form of $u_2$ now gives
\begin{eqnarray}
\label{positive side estimate}
\nonumber|\int_{\partial\Omega_{+,\epsilon_0}(\tilde\gamma)}\bar u_2\partial_\nu v| &=& |\int_{\partial\Omega_{+,\epsilon_0}(\tilde\gamma)}e^{-i\bar\zeta_2\cdot x}(e^{-i\bar\phi_2^\sharp}+ \bar r_2)\partial_\nu v|\\ &\leq& \| e^{-i\bar\phi_2^\sharp} +\bar r_2\|_{L^2(\partial\Omega)}\|e^{-s\tilde\gamma}\partial_\nu v\|_{L^2(\partial\Omega_{+,\epsilon_0}(\tilde\gamma))} 
\end{eqnarray}
Due to lemma \ref{L-infinity bound} and the standard theory of the restriction operator, the $\| e^{-i\bar\phi_2^\sharp} +\bar r_2\|_{L^2(\partial\Omega)}$ term in (\ref{positive side estimate}) satisfies the estimate
\[\| e^{-i\bar\phi_2^\sharp} +\bar r_2\|_{L^2(\partial\Omega)}\leq C + \|\bar r_2\|_{H^{\epsilon/4}(\partial\Omega)}\leq C(1+\|\chi_\Omega\bar r_2\|_{H^{1/2 + \epsilon/4}(\R^n)})\]
where $\chi_\Omega$ is a compactly supported smooth function that is identically one on $\Omega$. Since $\chi_\Omega$ is compactly supported, we have for all $-1<\delta<0$, 
\[\|\chi_\Omega\bar r_2\|_{H^{1/2 + \epsilon/4}(\R^n)}\leq C\|\chi_\Omega\bar r_2\|_{H_\delta^{1/2 + \epsilon/4}(\R^n)}\leq C\|\bar r_2\|_{H_\delta^{1/2 + \epsilon/4}(\R^n)}\] 
where the constant depends only on $\Omega$ and $\delta$. These facts combined with proposition \ref{existence of CGO} then gives the estimate
\begin{eqnarray}
\label{first multiplier}
\| e^{-i\bar\phi_2^\sharp} +\bar r_2\|_{L^2(\partial\Omega)}\leq C(1+ s^{1/2 - 3\epsilon/4})
\end{eqnarray}

Moving on to the $\|e^{-s\tilde\gamma}\partial_\nu v\|_{L^2(\partial\Omega_{+,\epsilon_0}(\tilde\gamma))}$ term in (\ref{positive side estimate}) we observe that proposition \ref{carleman} with Carleman weight $(-\tilde\gamma\cdot x)$ and operator $H_{W_2,q_2}$ applied to $v$ gives
\begin{eqnarray*}
\epsilon_0\|e^{-s\tilde\gamma}\partial_\nu v\|_{L^2(\partial\Omega_{+,\epsilon_0}(\tilde\gamma))}^2&\leq& (s^{-1}\|e^{-s\tilde\gamma\cdot x}H_{W_2,q_2}v\|_{L^2(\Omega)}^2 + \int_{\partial\Omega_{-,\epsilon_0}(\tilde\gamma)}|e^{-s\tilde\gamma\cdot x}\partial_\nu v|^2)\\
&\leq& s^{-1}\|e^{-s\tilde\gamma\cdot x} (W_1 - W_2)\cdot\nabla u_1\|_{L^2(\Omega)}^2 + \|\Lambda_{W_1,q_2}' - \Lambda_{W_2,q_2}'\|'^2 e^{2sD}\\\ && \ \ \ +s^{-1}\|e^{-s\tilde\gamma\cdot x}(\nabla\cdot(W_1 - W_2) + q_1 - q_2 + W_1^2 - W_2^2)u_1\|_{L^2(\Omega)}^2  \\
&\leq&  C(s+s^{-1} + \|\Lambda_{W_1,q_2}' - \Lambda_{W_2,q_2}'\|'^2 e^{2sD})
\end{eqnarray*}
for $s>>D$. Combining the above with (\ref{first multiplier}) and substitute the result into inequality (\ref{positive side estimate}) we get
\begin{eqnarray*}
|\int_{\partial\Omega_{+,\epsilon_0}(\tilde\gamma)}\bar u_2\partial_\nu v| 
\leq C(1 +s^{-1}+ s^{1-3\epsilon/4} + \|\Lambda_{W_1,q_2}' - \Lambda_{W_2,q_2}'\|'e^{sD}) 
\end{eqnarray*}
using the above inequality combined with (\ref{tail estimate}) and (\ref{the main term}) we can derive the following estimate from identity (\ref{fourier transform identity II})
\[|\int_{\Omega}\bar\mu\cdot(W_1 - W_2)e^{2i\xi\cdot x}|\leq C(\|\Lambda_{W_1,q_2}' - \Lambda_{W_2,q_2}'\|'e^{sD} +s^{-3\epsilon/4} + (\sqrt{1-\sqrt{1-|\xi|^2/s^2}})^{1-\frac{n}{p}})\]
Observe that since $Im(-\mu) = Im(\bar\mu)= (-\tilde\gamma)$, the same Carleman weight applies to give the same estimate for $-\mu$ in place of $\bar\mu$
\[|\int_{\Omega}-\mu\cdot(W_1 - W_2)e^{2i\xi\cdot x}|\leq C( \|\Lambda_{W_1,q_2}' - \Lambda_{W_2,q_2}'\|'e^{sD} +s^{-3\epsilon/4} + (\sqrt{1-\sqrt{1-|\xi|^2/s^2}})^{1-\frac{n}{p}})\]
Add the two estimates together and multiply by $|\xi_je_n -\xi_ne_j|$ we have the desired estimate.$\square$\\
We have established stability on the wedge $E_j$ for the Fourier transform of the component\\ $\partial_j(W_1 -W_2)_n - \partial_n(W_1 -W_2)_j$ of the curl. However we still have not established this estimate for the $(j,k)$ component when neither one of them is equal to $n$. But we will see in the next corollary that this follows immediately from lemma \ref{wedge estimate} - provided we take smaller sets than the ones we originally considered. We define the subset $\tilde E_j$ of $E_j$ by
\[\tilde E_j := \{\xi\in E_j \mid \xi_n\geq r_1/2 \xi_j\}\]
\begin{corollary}
\label{estimate for the other components}
The following estimate holds for all $\xi\in \tilde E_j\cap\tilde E_k$ such that $|\xi|\leq s$
\begin{eqnarray*}
&&|\int_\Omega e^{2i\xi\cdot x}\{\xi_j( W_1 - W_2)_k - \xi_k(W_1 - W_2)_j\}dx|\\
&&\leq  C|\xi|(e^{Ds}\|\Lambda'_{W_1,q_1} - \Lambda'_{W_1,q_2}\|' +(\sqrt{1-\sqrt{1-|\xi|^2/s^2}})^{1-\frac{n}{p}} + s^{-3\epsilon/4})
\end{eqnarray*}
\end{corollary}
\noindent{\bf Proof}\\
By lemma \ref{wedge estimate} we have for all $\xi$ in $\tilde E_j\cap\tilde E_k$
\begin{eqnarray*}
&&|\int_\Omega e^{2i\xi\cdot x}\{\xi_j( W_1 - W_2)_n - \xi_n(W_1 - W_2)_j\}dx|\\
&&\leq  C|\xi|(e^{Ds}\|\Lambda'_{W_1,q_1} - \Lambda'_{W_1,q_2}\|' +(\sqrt{1-\sqrt{1-|\xi|^2/s^2}})^{1-\frac{n}{p}} + s^{-3\epsilon/4})
\end{eqnarray*}
multiply this inequality by $|\frac{\xi_k}{\xi_n}|$ and use the definition of $\tilde E_k$ we have
\begin{eqnarray*}
&&|\int_\Omega e^{2i\xi\cdot x}\{\frac{\xi_j\xi_k}{\xi_n}( W_1 - W_2)_n - \xi_k(W_1 - W_2)_j\}dx|\\
&&\leq  \frac{C}{r_1}|\xi|(e^{Ds}\|\Lambda'_{W_1,q_1} - \Lambda'_{W_1,q_2}\|' +(\sqrt{1-\sqrt{1-|\xi|^2/s^2}})^{1-\frac{n}{p}} + s^{-3\epsilon/4})
\end{eqnarray*}
Switching the role of $j$ and $k$ we have the following estimate:
\begin{eqnarray*}
&&|\int_\Omega e^{2i\xi\cdot x}\{-\frac{\xi_j\xi_k}{\xi_n}( W_1 - W_2)_n + \xi_j(W_1 - W_2)_k\}dx|\\
&&\leq  \frac{C}{r_1}|\xi|(e^{Ds}\|\Lambda'_{W_1,q_1} - \Lambda'_{W_1,q_2}\|' +(\sqrt{1-\sqrt{1-|\xi|^2/s^2}})^{1-\frac{n}{p}} + s^{-3\epsilon/4})
\end{eqnarray*}
Add the two inequalities and we get the desired estimate.$\square$\\
\end{subsection}
\begin{subsection}{Proof of Theorem \ref{partial data stability for magnetic field}}
We first apply the result of Vessella to extend the estimate we have in a small wedge to an estimate on the entire ball of radius $R \leq s$. Define $f_R(\xi) := {\cal F}(\partial_j (\tilde W_1 - \tilde W_2)_k - \partial_k (\tilde W_1 - \tilde W_2)_j)(R\xi)$. Note that $f_R$ is analytic in $B_2$ and satisfies the estimate
\[|\partial_\xi^\alpha f_R(\xi)|\leq \frac{4Me^{nR}}{(D^{-1})^{|\alpha|}}\alpha!\ \ \ \ \forall \xi\in B_2\]
where $D := sup\{|x|\mid x\in\Omega\}$ and $M$ is the a-priori upper bound we have for the $W^{2,p}(\Omega)$ norm of $W_1$ and $W_2$. Proposition \ref{analytic cont} applies with $E = \tilde E_j \cap \tilde E_k \cap B_1$ to give the estimate
\[|f_R(\xi)| \leq (8e^{nR}M)^{1-\lambda}(\|f_R\|_{L^\infty(E)})^\lambda\ \ \ \forall \xi\in B_1\]
Here $\lambda\in(0,1)$ depends only on $r_1$, $\Omega$ and the dimension. By the fact that $\tilde E_j \cap \tilde E_k$ is a cone, we may apply corollary \ref{estimate for the other components} to obtain
\begin{eqnarray*}
\|f_R\|_{L^\infty(E)} &=& \|{\cal F}(\partial_j (\tilde W_1 - \tilde W_2)_k - \partial_k (\tilde W_1 - \tilde W_2)_j)\|_{L^\infty(B_R\cap\tilde E_j\cap\tilde E_k)}\\
&\leq& CR(e^{Ds}\|\Lambda'_{W_1,q_1} - \Lambda'_{W_2,q_2}\|' +(\sqrt{1-\sqrt{1-R^2/s^2}})^{1-\frac{n}{p}} + s^{-3\epsilon/4})
\end{eqnarray*}
as long as $s\geq R$. Combining the two inequalities and use the definition of $f_R$ we have
\begin{eqnarray}
\label{estimate on FT of magnetic field for partial data}
\|{\cal F}(\partial_j (\tilde W_1 - \tilde W_2)_k - \partial_k (\tilde W_1 - \tilde W_2)_j)\|_{L^\infty(B_R)}&\leq& CR^\lambda e^{nR(1-\lambda)}\times\\\nonumber &&(e^{Ds}\|\Lambda'_{W_1,q_1} - \Lambda'_{W_2,q_2}\|' +\sqrt{\frac{R^2}{s^2}}^{1-\frac{n}{p}} + s^{-3\epsilon/4})^\lambda
\end{eqnarray}
for all $R >0$ and $s \geq max\{R,1/h_0\}$. Here we used the fact that $R\leq s$ and for $t\in [0,1]$ we have the inequality $\sqrt{1-t}\geq 1-t$. With this estimate established for each component of ${\cal F}(d(\tilde W_1 - \tilde W_2))$, we now compute the $H^{-1}(\R^n)$ norm of $d(\tilde W_1 -\tilde W_2)$
\begin{eqnarray}
\label{H-1 norm of partial data magnetic field}
\nonumber\|d(\tilde W_1 &-& \tilde W_2)\|_{H^{-1}(\R^n)}^2 = \int_{B_R}\frac{|{\cal F}(d(\tilde W_1 - \tilde W_2))(\xi)|^2}{1+ |\xi|^2} d\xi +  \int_{\R^n\backslash B_R}\frac{|{\cal F}(d(\tilde W_1 - \tilde W_2))(\xi)|^2}{1+ |\xi|^2} d\xi\\
&\leq&C(R^{2\lambda+n}e^{2(1-\lambda)nR}(e^{2Ds}\|\Lambda'_{W_1,q_1} - \Lambda'_{W_2,q_2}\|'^2 +(R^2/s^2)^{\epsilon}+ s^{-3\epsilon/2})^\lambda + R^{-2})
\end{eqnarray}
here we have assumed without loss of generality that $\epsilon \leq 1-n/p$. We now choose $ s = R^2 R^{\frac{2\lambda + n}{\lambda\epsilon}}e^{2\frac{1-\lambda}{\epsilon\lambda}nR}$ so that
\[R^{2\lambda+n}e^{2(1-\lambda)nR}(e^{2Ds}\|\Lambda'_{W_1,q_1} - \Lambda'_{W_2,q_2}\|'^2 +(R^2/s^2)^{\epsilon}+ s^{-3\epsilon/2})^\lambda  \leq (e^{e^{KR}}\|\Lambda'_{W_1,q_1} - \Lambda'_{W_2,q_2}\|'^2 + e^{-(1-\lambda)R})^\lambda\]
for some $K>0$ depending only on $\Omega$, $\epsilon$, and $\lambda$. Substitute the above inequality into (\ref{H-1 norm of partial data magnetic field}) we get
\[\|d(\tilde W_1 - \tilde W_2)\|_{H^{-1}(\R^n)}^2 \leq C( (e^{e^{KR}}\|\Lambda'_{W_1,q_1} - \Lambda'_{W_2,q_2}\|'^2 + e^{-(1-\lambda)R})^\lambda + R^{-2})\]
and we are free to choose $R>0$. Assume without loss of generality that $\|\Lambda'_{W_1,q_1} - \Lambda'_{W_2,q_2}\|'\leq e^{-e}$ and set $R = \frac{1}{K}log(log\frac{1}{\|\Lambda'_{W_1,q_1} - \Lambda'_{W_2,q_2}\|'})$ the proof is complete.  
$\square$\\


\end{subsection}
\end{section}
\begin{section}{Stability for Electric Potential}
In this section we give a sketch of the proof of theorem \ref{partial data stability for electric potential}. We will omit the details since the proof is simply a combination of techniques employed in the proof of theorem \ref{partial data stability for magnetic field} and theorem \ref{potential estimate}. First we replace the magnetic potential $W_l$ by $W_l'$ as in the proof of the full data stability and observe that $\Lambda_{W_l',q_l}' = \Lambda_{W_l,q_l}'$. Define $u_1$, $u_2$ and $v$ as in proof of lemma \ref{wedge estimate}. Following the same steps as in proof of theorem \ref{potential estimate}, we have that for all $\xi\in \bigcup_{j = 1}^{n-1}E_j$ such that $|\xi|\leq R$
\begin{eqnarray}
\label{partial data potential FT}
\nonumber|{\cal F}(q_1 - q_2)(2\xi)|&\leq& C(e^{2Ds}\|\Lambda_{W_1,q_1}'-\Lambda_{W_2,q_2}'\|'  + \|d(W_1-W_2)\|_{H^1(\Omega)^*}^{\lambda/2}s + s^{-\epsilon}\\
&&\ \ \ \ \ \ \ \ \ \ \ \ \ \ \ \ \ \ \ \ \ \ \ \ \ \ \ \ \ \ \ \ \ \ \ \ \ \ \ \ \ \ \ \ \ \ \ +\int_{\partial\Omega_{+,\epsilon_0}(\tilde\gamma)(\tilde\gamma)}|\bar u_2\partial_\nu v|)
\end{eqnarray} 
As in the proof of lemma \ref{wedge estimate}
\begin{eqnarray}
\label{intermediate term}
\int_{\partial\Omega_{+,\epsilon_0}(\tilde\gamma)}|\bar u_2\partial_\nu v|\leq C(1 + s^{1-3\epsilon/4})\|e^{-s\tilde\gamma}\partial_\nu v\|_{L^2(\partial\Omega_{+,\epsilon_0}(\tilde\gamma))} 
\end{eqnarray}

Apply Carleman estimate to the second term on the RHS of (\ref{intermediate term}) we have

\begin{eqnarray}
\label{carleman II}
\nonumber\epsilon_0\|e^{-s\tilde\gamma}\partial_\nu v\|_{L^2(\partial\Omega_{+,\epsilon_0}(\tilde\gamma))}^2 &\leq&  s^{-1}(\|e^{-s\tilde\gamma\cdot x}(W_1' - W_2')\cdot\nabla u_1\|_{L^2(\Omega)}^2+C) + \|\Lambda_{W_1,q_2}' - \Lambda_{W_2,q_2}'\|'^2 e^{2sD}\\\nonumber
&\leq& s^{-1}(\|W_1' - W_2'\|_{L^\infty(\Omega)}s^2+C) + \|\Lambda_{W_1,q_2}' - \Lambda_{W_2,q_2}'\|'^2 e^{2sD}\\\nonumber
&\leq& s^{-1}(\|W_1' - W_2'\|_{L^\infty(\Omega)}s^2+C) + \|\Lambda_{W_1,q_2}' - \Lambda_{W_2,q_2}'\|'^2 e^{2sD}\\
&\leq& s^{-1}C(\|d(W_1 - W_2)\|_{H^1(\Omega)^*}s^2+1) + \|\Lambda_{W_1,q_2}' - \Lambda_{W_2,q_2}'\|'^2 e^{2sD}
\end{eqnarray}
where the last inequality comes from applying lemma \ref{bounded inverse of d} to $W_1' - W_2'$ then use interpolation inequalities. Combining (\ref{carleman II}), (\ref{intermediate term}), and (\ref{partial data potential FT}) we have the following estimate for the Fourier transform of $(q_1 - q_2)$ in the wedge $\xi\in\bigcup_{j = 1}^{n-1}E_j$
\begin{eqnarray*}
|{\cal F}(q_1 - q_2)(2\xi)|\leq C(e^{2Ds}\|\Lambda_{W_1,q_1}'-\Lambda_{W_2,q_2}'\|' + s^{-3\epsilon/4} +\|d(W_1-W_2)\|_{H^1(\Omega)^*}^{\lambda/2}s^2)
\end{eqnarray*}
Theorem \ref{partial data stability for magnetic field} applies to give
\begin{eqnarray*}
|{\cal F}(q_1 - q_2)(2\xi)|&\leq& C(e^{2Ds}\|\Lambda_{W_1,q_1}'-\Lambda_{W_2,q_2}'\|' + s^{-3\epsilon/4} \\&&+(\|\Lambda_{W_1,q_1}'-\Lambda_{W_2,q_2}'\|'^\lambda + \frac{1}{|log( |log \|\Lambda_{W_1,q_1}'-\Lambda_{W_2,q_2}'\|'|)|})^{\lambda/2}s^2)
\end{eqnarray*}
for all $\xi \in \bigcup_{j = 1}^{n-1}E_j$. The theorem now follows from proposition \ref{analytic cont} and standard computation.$\square$

\end{section}
\end{section}
  
\end{document}